\documentclass{IJFS}

\usepackage[utf8]{inputenc}
\usepackage{amsfonts,amsthm,amsmath,amssymb,dsfont}
\usepackage{amsmath,    enumitem}
\usepackage{graphicx}
\usepackage[all]{xy}
\usepackage{subfigure}
\usepackage{tikz}
\usetikzlibrary[arrows,shapes,graphs,decorations.pathmorphing,backgrounds,positioning,fit,petri]
\usepackage{amscd, eufrak}

 \usepackage{amssymb}

 \usepackage{mathrsfs}
 \usepackage{amsfonts}

 \usepackage{verbatim}

\usepackage{color,colortbl}

\definecolor{OrcidGreen}{rgb}{0.6,0.8,0.2}

 \usepackage[singlelinecheck=false]{caption}
 \usepackage[linesnumbered,ruled]{algorithm2e}
 \usepackage[colorlinks=true,linkcolor=black,citecolor=black,urlcolor=blue]{hyperref}

\usepackage[square,numbers,sort&compress]{natbib}
\usepackage{url}
\usepackage{hyphenat}
\usepackage{parskip}

\usepackage{tabularray}

\usepackage{upgreek}
\usepackage{breqn}
\usepackage{mathtools}

\usepackage{pifont}
\usepackage{fancyhdr}
\usepackage{subfiles}
\usepackage{fontawesome5}

\usepackage{amsthm}

\newtheorem{theorem}{Theorem}[section]
\newtheorem{lemma}[theorem]{Lemma}

\newtheorem{corollary}[theorem]{Corollary}

\newtheorem{definition}[theorem]{Definition}
\newtheorem{example}[theorem]{Example}
\newtheorem{remark}[theorem]{Remark}
\newtheorem{therule}[theorem]{Rule}

\begin{document}

\addtocounter{page}{0}

\titleijfs{On the nonlinear programming problems subject to a system of generalized bipolar fuzzy relational equalities defined with continuous t-norms \vspace{-2mm}}
\author[1 \href{https://orcid.org/0000-0002-9224-8470}{\textcolor{OrcidGreen}{\faOrcid}}]{A. Ghodousian}
\author[2 \href{https://orcid.org/0009-0006-4046-2674}{\textcolor{OrcidGreen}{\faOrcid}}]{M. S. Chopannavaz}
\affil[1]{Faculty of Engineering Science, College of Engineering, University of Tehran, P.O.Box 11365-4563, Tehran, Iran.}
\affil[2]{Department of Engineering Science, College of Engineering, University of Tehran, Tehran, Iran.}
\emails{a.ghodousian@ut.ac.ir, chopannavaz@ut.ac.ir}
\CorrespondAuthor{A. Ghodousian}
\oddPageHead{On the nonlinear programming problems subject to a system of generalized BFRE defined with continuous t-norms}
\evenPageHead{A. Ghodousian, M. S. Chopannavaz}

\abstractijfs{
As a starting point, this paper develops the system of bipolar fuzzy relational equations (FRE) to the most general case, where bipolar FREs are defined by an arbitrary continuous t-norm. Due to the fact that fuzzy relational equations are special cases of bipolar FREs, the proposed system can also be viewed as a generalization of traditional FREs, in which the fuzzy composition can be defined by a continuous t-norm. In order to determine the feasibility of the proposed system, some necessary and sufficient conditions are presented for studying continuous bipolar FREs. This is followed by a complete analysis of the set of feasible solutions to the problem. Contrary to FREs and bipolar FREs defined by continuous Archimedean t-norms, the feasible solutions set of generalized bipolar FREs consists of a finite number of compact sets that are not necessarily connected. Further, five techniques have been outlined in an attempt to simplify the current problem, and then an algorithm has been presented to find the feasible region of the problem. Next, we present a class of optimization models subject to continuous bipolar FRE constraints, in which the objective function incorporates a wide range of (non)linear functions, such as maximum functions, geometric mean functions, log-sum-exp functions, maximum eigenvalues of symmetric matrices, support functions for sets, etc. Considering that the problem has a finite number of local optimal solutions, the global optimal solution can always be obtained by choosing the point with the minimum objective value among these local optimal solutions. Lastly, as a means to illustrate the definitions, theorems, and algorithms presented in the paper, a step-by-step example is presented in several sections, in which the constraints are a system of bipolar FREs defined by the Dubois-Prade t-norm, which is a continuous non-Archimedean t-norm.
}

\keywordsijfs{Bipolar fuzzy relational equations, continuous t-norms, Max-continuous compositions, Global optimization, Non-linear optimization.
}

\Vol{22}	\No{3}	\Year{2025}	\Pages{21-38}
\Received{December 2024}	\Revised{February 2025}		\Accepted{May 2025}
\DOI{https://doi.org/10.22111/ijfs.2025.50433.8905}

{\let\newpage\relax\vspace*{.1mm}
\noindent\rule{\textwidth}{.3mm}\maketitle}

\thispagestyle{fancylogo}

\mabstract

% ==================================== section1_Introduction ====================================
\section{Introduction}\label{sec_intro}
The theory of fuzzy relational equations (FRE) as a generalized version of Boolean relation equations was firstly proposed by Sanchez and was applied to problems related to the medical diagnosis \cite{ref_43}. Pedrycz categorized and extended two ways of the generalizations of FRE in terms of sets under discussion and various operations which are taken account \cite{39}. Since then, FRE was applied in many other fields such as fuzzy control, prediction of fuzzy systems, fuzzy decision making, fuzzy pattern recognition, image compression and reconstruction, fuzzy clustering and so on. Generally, when rules of inference are applied and their corresponding consequences are known, the problem of determining antecedents is simplified and mathematically reduced to solving an FRE \cite{37}. Nowadays, it is well known that many of the issues associated to the body knowledge can be treated as FRE problems \cite{38}. Because of important applications in various practical fields, many scholars have focused on the theoretical research of FRE, including its resolution approach and some specific optimization problems with FRE constraints. 

The solvability identification and finding set of solutions are the primary, and the most fundamental, matter concerning the FRE problems. Di Nola et al. proved that the solution set of FRE (if it is non-empty), defined by continuous max-t-norm composition is often a non-convex set. This non-convex set is completely determined by one maximum solution and a finite number of minimal solutions \cite{8}. Such non-convexity property is one of two bottlenecks making a major contribution towards an increase in complexity of FRE-related problems, particularly, in the optimization problems subjected to a system of fuzzy relations. Another bottleneck point is concerned with detecting the minimal solutions for FREs. Chen and Wang \cite{4} presented an algorithm for obtaining the logical representation of all minimal solutions and deduced that a polynomial-time algorithm with the ability to find all minimal solutions of FRE (with max-min composition) may not exist. Also, Markovskii showed that solving max-product FRE is closely related to the covering problem which is a type of NP-hard problem \cite{36}. In fact, the same result holds true for a more general t-norms instead of the minimum and product operators \cite{5,32,33}. Over the past decades, the solvability of FRE which is defined using different max-t compositions have been investigated by many researchers \cite{20,21,23,41,45,47,48,51,53}. Moreover, some other researchers have worked on introducing novel concept and at times improving some of the existing theoretical aspects and applications of fuzzy relational inequalities (FRI) \cite{18,22,23,24,29,55}. For example, Li and Yang \cite{29} studied FRI with addition-min composition and presented an algorithm to search for minimal solutions. They applied FRI to data transmission mechanism in a BitTorrent-like Peer-to-Peer file sharing systems. In \cite{18}, the authors focused on the study of a mixed fuzzy system formed by two fuzzy relational inequalities $A \varphi x \leq b^{1}$ and $D \varphi x \geq b^{2}$, where $\varphi$ is an operator with (closed) convex solutions.

The optimization problem subject to FRE and FRI is one of the most interesting and on-going research topics amongst similar problems \cite{3,13,17,18,20,21,22,23,24,25,26,28,34,42,46,49,55}. Many methods were designed based on the translation of the main problem into an integer linear programming problem which is then solved using well-developed techniques. On the contrary, other algorithms benefit the resolution of the feasible region, some necessary and sufficient conditions for the optimality and simplification processes. The most methods of this category are based on analytical results provided mainly by Sanchez \cite{44} and Pedrycz \cite{40}. For instance, Fang and Li converted a linear optimization problem subjected to FRE constraints with max-min operation into an integer programming problem and solved it by a branch-and-bound method using jumptracking technique \cite{14}. Wu et al. worked on improvement of the method employed by Fang and Li this was done by decreasing the search domain and presented a simplification process by three rules which resulted from a necessary condition \cite{50}. Chang and Shieh presented new theoretical results concerning the linear optimization problem constrained by fuzzy max-min relation equations \cite{3}. They improved an upper bound on the optimal objective value, some rules for simplifying the problem and proposed a rule for reducing the solution tree. In \cite{27}, an application of optimizing the linear objective with max-min composition was employed for the streaming media provider seeking a minimum cost while fulfilling the requirements assumed by a three-tier framework. Linear optimization problem was further investigated by numerous scholars focusing on max-product operation \cite{26,35}. Loetamonphong and Fang defined two subproblems by separating negative from non-negative coefficients in the objective function, and then obtained an optimal solution by combining the optimal solutions of the two sub-problems \cite{35}. Moreover, generalizations of the linear optimization problem with respect to FRE have been studied; this was done through replacement of max-min and max-product compositions with different fuzzy compositions such as max-average composition \cite{49} or max-t-norm composition \cite{20,21,23,25,28,46}. For example, Li and Fang solved the linear optimization problem subjected to a system of sup-t equations by reducing it to a 0-1 integer optimization problem \cite{28}. In \cite{25}, a method was presented for solving linear optimization problems with the max-Archimedean tnorm fuzzy relation equation constraint. In \cite{46}, the authors solved the same problem whit continuous Archimedean t-norm, and to obtain some optimal variables, they used the covering problem rather than the branch-and-bound methods.

Recently, many interesting forms of generalizations of the linear programming applied to the system of fuzzy relations have been introduced, and developed based on composite operations used in FRE, fuzzy relations used in the definition of the constraints, some developments on the objective function of the problems and other ideas \cite{7,9,15,23,31,34,52,56}. For example, Wu et al. represented an efficient method to optimize a linear fractional programming problem under FRE with max-Archimedean t-norm composition \cite{52}. Dempe and Ruziyeva generalized the fuzzy linear optimization problem by considering fuzzy coefficients \cite{7}. In addition, Dubey et al. studied linear programming problems involving interval uncertainty modeled using intuitionistic fuzzy set \cite{9}. Yang \cite{56} studied the optimal solution of minimizing a linear objective function subject to a FRI where the constraints defined as $\sum_{j=1}^{n} \min \left\{a_{i j}, x_{j}\right\} \geq b_{i}$ for $i=1, \ldots, m$. Also, in \cite{55}, the authors introduced the latticized linear programming problem subject to max-product fuzzy relation inequalities with application in the optimization management model for wireless communication emission base stations. The latticized linear programming problem was defined by minimizing the objective function $z(x)=\min_{j=1}^{n}\left\{x_{j}\right\}$ subject to the feasible region $X(A, b)=\left\{x \in[0,1]^{n}: A \circ x \geq b\right\}$ where "$\circ$" denotes fuzzy max-product composition. They also presented an algorithm based on the resolution of the feasible region.

The concept of FRE was also generalized to a so-called bipolar fuzzy relation equation. Bipolarity exists widely in human understanding of information and preference \cite{11}. Dubois and Prade provided an overview of the asymmetric bipolar representation of positive and negative information in possibility theory \cite{12}. They showed that the possibility theory framework is convenient for handling bipolar representations, that was applied to distinguish between negative and positive information in preference modeling \cite{11,12}. The linear optimization of bipolar FRE was also the focus of study carried out by some researchers where FRE was defined with max-min (with application in public awareness of the products for a supplier) \cite{15,30}, maxproduct \cite{6} and max-Lukasiewicz composition \cite{31,34,54,57}. In \cite{15}, the concept of bipolar FRE was firstly proposed with max-min composition where the constraints are expressed as $\max _{j=1}^{n}\left\{\max \left\{\min \left\{a_{i j}^{+}, x_{j}\right\}, \min \left\{a_{i j}^{-}, 1-x_{j}\right\}\right\}\right\}$ for $i=1, \ldots, m$, where $a_{i j}^{+}, a_{i j}^{-}$, $x_{j} \in[0,1]$. Similarly, in \cite{31}, the authors introduced a linear optimization problem subjected to a system of bipolar FRE defined as $X\left(A^{+}, A^{-}, b\right)=\left\{x \in[0,1]^{m}: x \circ A^{+} \vee \tilde{x} \circ A^{-}=b\right\}$ where $\tilde{x}_{i}=1-x_{i}$ for each component of $\tilde{x}=\left(\tilde{x}_{i}\right)_{1 \times m}$ and the notations "$\vee$" and "$\circ$" denote max operation and the max-Lukasiewicz composition, respectively. They translated the original problem into a $0-1$ integer linear problem. In a separate, the foregoing bipolar linear optimization problem was solved by an analytical method based on the resolution and some structural properties of the feasible region (using a necessary condition for characterizing an optimal solution and a simplification process for reducing the problem) \cite{34}. However, resolution method for obtaining the complete solution set of bipolar max-Lukasiewicz FRE was not found in the mentioned works \cite{54}. Yang \cite{54} studied bipolar max-Lukasiewicz FRE and showed that the complete solution set of the system is fully determined by finite conservative bipolar paths. Also, in \cite{16}, the authors studied bipolar FREs defined by max-strict t-norms. They categorized the constraints of the problem into four groups, each requiring a different approach to determine the corresponding feasible solution set.

In this paper, a wide class of bipolar FRE linear optimization problems is studied as the following mathematical model:
\begin{equation}\label{eq_1}
	\begin{array}{ll}
		\min & f(x) \\
		& A^{+} \varphi x \vee A^{-} \varphi(\mathbf{1}-x) = b \\
		& x \in[0,1]^{n}
	\end{array}
\end{equation}
where $A^{+}=\left(a_{i j}^{+}\right)_{m \times n}$ and $A^{-}=\left(a_{i j}^{-}\right)_{m \times n}$ are fuzzy matrices and $b=\left(b_{i}\right)_{m \times 1}$ is a fuzzy vector such that $0 \leq a_{i j}^{+} \leq 1$, $0 \leq a_{i j}^{-} \leq 1$ and $0 \leq b_{i} \leq 1$ for each $i \in \mathscr{I}=\{1,2, \ldots, m\}$ and each $j \in \mathscr{J}=\{1,2, \ldots, n\}$, respectively. Also, the sum vector 1 (a vector with each component equal to one) and the unknown vector $x=\left(x_{j}\right)_{n \times 1}$ are in $\mathbb{R}^{n}$. Moreover, $f: \mathbb{R}^{n} \rightarrow \mathbb{R}$ is a continuous function that is non-decreasing (non-increasing) in $x_{j}$, $\forall j \in \mathscr{J}^{+}\left(\forall j \in \mathscr{J}^{-}\right)$, where $\mathscr{J}^{+} \subseteq \mathscr{J}$, $\mathscr{J}^{-} \subseteq \mathscr{J}$, $\mathscr{J}^{+} \cap \mathscr{J}^{-}=\varnothing$ and $\mathscr{J}^{+} \cup \mathscr{J}^{-}=\mathscr{J}$. Additionally, $\varphi$ denotes an arbitrary continuous t-norm and the constraints in Problem (\ref{eq_1}) mean:
\begin{equation}\label{eq-2}
	\max _{j=1}^{n}\left\{\max \left\{\varphi\left(a_{i j}^{+}, x_{j}\right), \varphi\left(a_{i j}^{-}, 1-x_{j}\right)\right\}\right\}=b_{i} \quad, \forall i \in \mathscr{I}
\end{equation}
In the current paper, it is shown that many important functions (from the theoretical or practical point of view) satisfy the assumptions expressed for the objective function of Problem (\ref{eq_1}); for example, we can refer to linear function, maximum function, geometric mean function, log-sumexp function, maximum eigenvalue of a symmetric matrix, support function of a set, sum of $r$ largest components, sum of logarithms, any norm function and norm of matrices, etc. On the other hand, it is shown that the feasible region of Problem (\ref{eq_1}), unlike FREs \cite{18} and bipolar FREs defined by Archimedean t-norms \cite{16}, is determined as the union of a finite number of compact sets that are not necessarily connected. However, it is proved that problem (\ref{eq_1}) has a limited number of local optimal solutions and the global optimum is one of the local optima with the smallest objective function value.

The rest of the paper is organized as follows. Some preliminary definitions, concepts and properties of continuous bipolar FREs are presented in Section \ref{sec-2}. In Section \ref{sec-3}, the feasible solution set of the original problem is characterized and two necessary conditions are derived to identify the feasibility of the problem. Moreover, a necessary and sufficient condition is given to guarantee that a given point is feasible to Problem (\ref{eq_1}) or not. Section \ref{sec-4} describes five simplification rules introduced to accelerate the resolution process by reducing the size of the problem, followed by an algorithm for resolution of the feasible region of the current optimization problems. In Section \ref{sec-5}, Problem (\ref{eq_1}) is solved by finding the local optima and the global optimal solution as the local optimum with the smallest objective function value.

% ==================================== section2_Preliminary definitions and properties ====================================
\section{Preliminary definitions and properties}\label{sec-2}
This section initially describes the feasible solution set of the equation $\varphi(a, x)=b$ where $a$ and $b$ are two fixed scalars in $[0,1]$, $x \in[0,1]$ and $\varphi$ is an arbitrary continuous t-norm. Then, the feasible solution set of the equation $\max \left\{\varphi\left(a_{i j}^{+}, x_{j}\right), \varphi\left(a_{i j}^{-}, 1-x_{j}\right)\right\}=b_{i}$ is completely characterized for each $i \in \mathscr{I}$ and each $j \in \mathscr{J}$. Subsequently, according to the results obtained in this section, the feasible region of Problem (\ref{eq_1}) is determined in the next section. For the sake of simplicity, let $S_{i}$ denote the feasible solution set of the $i$'th equation; that is, $S_{i}=\left\{x \in[0,1]^{n}: \max _{j=1}^{n}\left\{\max \left\{\varphi\left(a_{i j}^{+}, x_{j}\right), \varphi\left(a_{i j}^{-}, 1-x_{j}\right)\right\}\right\}=b_{i}\right\}$. Also, let $S\left(A^{+}, A^{-}, b\right)$ denote the feasible solution set for Problem (\ref{eq_1}). So, it is clear that $S\left(A^{+}, A^{-}, b\right)=\bigcap_{i \in \mathscr{J}} S_{i}$.

\begin{definition}\label{def-1}
	Suppose that $a$ and $b$ are two fixed scalers in $[0,1]$. We define $S=\{x \in[0,1]: \varphi(a, x)=b\}$ and $I=\{x \in[0,1]: \varphi(a, x) \leq b\}$. Also, if $S \neq \varnothing$, define $l=\inf S$ and $u=\sup S$. 
\end{definition}

It should be noted that if $S \neq \varnothing$ in Definition \ref{def-1}, then both $l=\inf S$ and $u=\sup S$ exist and we also have the following results \cite{18}:

\begin{lemma}\label{lm-1}
	Suppose that $a$ and $b$ are two fixed scalers in [0,1] and $\varphi$ is a continuous t-norm. Then,
	\begin{equation}\label{eq_3}
		S=\left\{\begin{array}{cc}
			{[l, u]} & , a \geq b  \\
			\varnothing & , a<b
		\end{array}, \quad I= \begin{cases}{[0, u]} & , a \geq b \\
			{[0,1]} & , a<b\end{cases}\right.
	\end{equation}
\end{lemma}

For the special cases of Lemma \ref{lm-1}, we refer the reader to \cite{16} where $\varphi$ is a continuous strict tnorm. Table \ref{tbl_a1} (see Appendix \ref{sec_appA}) presents some commonly used continuous t-norms. For each t-norm in Table \ref{tbl_a1}, values $l$ and $u$ discussed in Lemma \ref{lm-1} is presented in Table \ref{tbl_a2}.

\begin{definition}\label{def-2}
	For each $i \in \mathscr{J}$ and each $j \in \mathscr{J}$, we define $S_{i j}^{+}=\{x \in[0,1]: \varphi(a_{i j}^{+}, x)=b_{i}\}$ and $S_{i j}^{-}=\{x \in[0,1]: \varphi(a_{i j}^{-}, 1-x)=b_{i}\}$; that is, $S_{i j}^{+}$ and $S_{i j}^{-}$ denote the feasible solution sets of the equations $\varphi(a_{i j}^{+}, x)=b_{i}$ and $\varphi(a_{i j}^{-}, 1-x)=b_{i}$, respectively. Furthermore, define $I_{i j}^{+}=\{x \in[0,1]: \varphi(a_{i j}^{+}, x) \leq b_{i}\}$ and $I_{i j}^{-}=\{x \in[0,1]: \varphi(a_{i j}^{-}, 1-x) \leq b_{i}\}$.
\end{definition}

According to Definition \ref{def-2}, it is easy to verify that $S_{i j}^{+} \subseteq I_{i j}^{+}$and $S_{i j}^{-} \subseteq I_{i j}^{-}, \forall i \in \mathscr{I}$ and $\forall j \in \mathscr{J}$. Also, the following results are directly obtained from (\ref{eq_3}) and Definition \ref{def-2}:
\begin{equation}\label{eq-4}
	\begin{alignedat}{2}
		S_{i j}^{+}&=\begin{cases}
			{\left[l_{i j}^{+}, u_{i j}^{+}\right]} & , a_{i j}^{+} \geq b_{i} \\
			\varnothing & , a_{i j}^{+}<b_{i}
		\end{cases} \qquad
		&\quad ,I_{i j}^{+}&=\begin{cases}{\left[0, u_{i j}^{+}\right]} & , a_{i j}^{+} \geq b_{i} \\ 
			{[0,1]} & , a_{i j}^{+}<b_{i}\end{cases}\\
		S_{i j}^{-}&=\begin{cases}
			{\left[1-u_{i j}^{-}, 1-l_{i j}^{-}\right]} & , a_{i j}^{-} \geq b_{i} \\
			\varnothing & , a_{i j}^{-}<b_{i}
		\end{cases}
		&\quad ,I_{i j}^{-}&=\begin{cases}{\left[1-u_{i j}^{-}, 1\right]} & , a_{i j}^{-} \geq b_{i} \\
			{[0,1]} & , a_{i j}^{-}<b_{i}\end{cases}
	\end{alignedat}
\end{equation}
where $l_{i j}^{+}=\inf S_{i j}^{+}$, $u_{i j}^{+}=\sup S_{i j}^{+}$, $l_{i j}^{-}=\inf \left\{x \in[0,1]: \varphi\left(a_{i j}^{-}, x\right)=b_{i}\right\}$ and $u_{i j}^{-}=\sup \left\{x \in[0,1]: \varphi\left(a_{i j}^{-}, x\right)=b_{i}\right\}$.

\begin{definition}\label{def-3}
	For each $i \in \mathscr{I}$ and each $j \in \mathscr{J}$, we define $S_{i j}=\left\{x \in[0,1]: \max \left\{\varphi\left(a_{i j}^{+}, x\right), \varphi\left(a_{i j}^{-}, 1-x\right)\right\}=b_{i}\right\}$ and
	$I_{i j}=\left\{x \in[0,1]: \max \left\{\varphi\left(a_{i j}^{+}, x\right), \varphi\left(a_{i j}^{-}, 1-x\right)\right\} \leq b_{i}\right\}$.
\end{definition}

\begin{lemma}\label{lm-2}
	\textbf{(a)} $I_{i j}=I_{i j}^{+} \cap I_{i j}^{-}$. \textbf{(b)} $S_{i j}=I_{i j} \cap\left(S_{i j}^{+} \cup S_{i j}^{-}\right)$.
\end{lemma}

\begin{proof}
	\textbf{(a)} The proof is easily resulted from Definitions \ref{def-2} and \ref{def-3}. \textbf{(b)} According to Definition \ref{def-3}, $x \in S_{i j}$ iff $\varphi\left(a_{i j}^{+}, x_{j}\right) \leq b_{i}$, $\varphi\left(a_{i j}^{-}, 1-x_{j}\right) \leq b_{i}$ (i.e., $x \in I_{i j}^{+} \cap I_{i j}^{-}=I_{i j}$) and at least one of the two equalities $\varphi\left(a_{i j}^{+}, x_{j}\right)=b_{i}$ (i.e., $x \in S_{i j}^{+}$) and $\left.\varphi\left(a_{i j}^{-}, 1-x_{j}\right)\right\}=b_{i}$ (i.e., $x \in S_{i j}^{-}$) holds.
\end{proof}

Based on (\ref{eq-4}) and Lemma \ref{lm-2}, we can obtain the following results, which characterize sets $S_{i j}$ and $I_{i j}$ for all cases.

\begin{corollary}\label{corl-1}
	Suppose that $i \in \mathscr{I}, j \in \mathscr{J}$ and $\varphi$ is a continuous t-norm.
	\begin{itemize}
		\item[\textbf{(a)}] If $a_{i j}^{+}<b_{i}$ and $a_{i j}^{-}<b_{i}$, then $I_{i j}=[0,1]$ and $S_{i j}=\varnothing$.
		
		\item[\textbf{(b)}] If $a_{i j}^{+} \geq b_{i}$ and $a_{i j}^{-}<b_{i}$, then $I_{i j}=\left[0, u_{i j}^{+}\right]$and $S_{i j}=\left[l_{i j}^{+}, u_{i j}^{+}\right]$.
		
		\item[\textbf{(c)}] If $a_{i j}^{-} \geq b_{i}$ and $a_{i j}^{+}<b_{i}$, then $I_{i j}=\left[1-u_{i j}^{-}, 1\right]$ and $S_{i j}=\left[1-u_{i j}^{-}, 1-l_{i j}^{-}\right]$.
		
		\item[\textbf{(d)}] If $a_{i j}^{+} \geq b_{i}$ and $a_{i j}^{-} \geq b_{i}$, then $I_{i j}=\left[1-u_{i j}^{-}, u_{i j}^{+}\right]$ and
		$$
		S_{i j}= \begin{cases}{\left[1-u_{i j}^{-}, 1-l_{i 			j}^{-}\right] \bigcup\left[l_{i j}^{+}, u_{i j}^{+}\right]} & , l_{i j}^{+}>1-l_{i j}^{-} \\ {\left[1-u_{i j}^{-}, u_{i j}^{+}\right]} & , l_{i j}^{+} \leq 1-l_{i j}^{-}\end{cases}.
		$$
	\end{itemize}
\end{corollary}

\begin{definition}\label{def-4}
	For each $j \in \mathscr{J}$, we define $I^{-}(j)=\left\{i \in \mathscr{I}: a_{i j}^{-} \geq b_{i}\right\}$, $I^{+}(j)=\left\{i \in \mathscr{I}: a_{i j}^{+} \geq b_{i}\right\}$ and $I_{j}=\bigcap_{i \in \mathscr{I}} I_{i j}$. Also, define $S_{i j}^{\prime}=S_{i j} \cap I_{j}$, $\forall i \in \mathscr{I}$ and $\forall j \in \mathscr{J}$.
\end{definition}

\begin{remark}\label{rmk-1}
	According to Corollary \ref{corl-1} and Definition \ref{def-4}, it is concluded that $I_{j}=\left[L_{j}, U_{j}\right]$, $\forall j \in \mathscr{J}$, where
	$$
	L_{j}=\left\{\begin{array}{ll}
		\max _{i \in I^{-}(j)}\left\{1-u_{i j}^{-}\right\} & , I^{-}(j) \neq \varnothing \\
		0 & , I^{-}(j)=\varnothing
	\end{array} \quad, \quad U_{j}= \begin{cases}\min _{i \in I^{+}(j)}\left\{u_{i j}^{+}\right\} & , I^{+}(j) \neq \varnothing \\
		1 & , I^{+}(j)=\varnothing\end{cases}\right.
	$$
	Moreover, for each $i \in \mathscr{I}$ and $j \in \mathscr{J}$ such that $S_{i j}^{\prime} \neq \varnothing$, we have
	$$
	S_{i j}^{\prime}= \begin{cases}{\left[\max \left\{L_{j}, l_{i j}^{+}\right\}, \min \left\{U_{j}, u_{i j}^{+}\right\}\right]} & , a_{i j}^{+} \geq b_{i}, ~a_{i j}^{-}<b_{i} \\ {\left[\max \left\{L_{j}, 1-u_{i j}^{-}\right\}, \min \left\{U_{j}, 1-l_{i j}^{-}\right\}\right]} & , a_{i j}^{+}<b_{i}, ~a_{i j}^{-} \geq b_{i} \\ {\left[\max \left\{L_{j}, 1-u_{i j}^{-}\right\}, \min \left\{U_{j}, 1-l_{i j}^{-}\right\}\right] \cup\left[\max \left\{L_{j}, l_{i j}^{+}\right\}, \min \left\{U_{j}, u_{i j}^{+}\right\}\right]} & , a_{i j}^{+} \geq b_{i}, ~a_{i j}^{-} \geq b_{i}, ~l_{i j}^{+}>1-l_{i j}^{-} \\ {\left[\max \left\{L_{j}, 1-u_{i j}^{-}\right\}, \min \left\{U_{j}, u_{i j}^{+}\right\}\right]} & , a_{i j}^{+} \geq b_{i}, ~a_{i j}^{-} \geq b_{i}, ~l_{i j}^{+} \leq 1-l_{i j}^{-}\end{cases}
	$$
\end{remark}

\begin{example}\label{ex-1}
	Consider Problem (\ref{eq_1}) with the Dubois-Prade operator where
	$$A^{+}=\left[\begin{array}{ccccccccc}0.54 & 0.48 & 0.80 & 0.63 & 0.70 & 0.35 & 0.56 & 0.29 & 0.69 \\ 0.20 & 0.06 & 0.01 & 0.03 & 0.00 & 0.04 & 0.50 & 0.00 & 0.09 \\ 0.72 & 0.23 & 0.75 & 0.44 & 0.38 & 0.61 & 0.51 & 0.80 & 0.67 \\ 0.83 & 1.00 & 0.30 & 0.90 & 0.89 & 0.79 & 0.62 & 0.41 & 0.86 \\ 0.13 & 0.10 & 0.00 & 0.15 & 0.11 & 0.04 & 0.00 & 0.07 & 0.19 \\ 0.28 & 0.43 & 0.35 & 0.28 & 0.40 & 0.22 & 0.18 & 0.50 & 0.00 \\ 0.33 & 0.60 & 0.54 & 0.58 & 0.14 & 0.80 & 0.49 & 0.26 & 0.39\end{array}\right]$$ $$A^{-}=\left[\begin{array}{ccccccccc}0.65 & 0.51 & 0.70 & 0.26 & 0.90 & 0.46 & 0.68 & 0.16 & 0.29 \\ 0.10 & 0.20 & 0.00 & 0.06 & 0.03 & 0.00 & 0.05 & 0.00 & 0.00 \\ 0.13 & 0.63 & 0.74 & 0.25 & 0.66 & 0.73 & 0.39 & 0.80 & 0.90 \\ 0.81 & 0.80 & 0.92 & 0.90 & 0.78 & 0.88 & 0.95 & 0.57 & 0.18 \\ 0.17 & 0.25 & 0.09 & 0.18 & 0.40 & 0.00 & 0.19 & 0.08 & 0.00 \\ 0.00 & 0.29 & 0.33 & 0.47 & 0.27 & 0.34 & 0.15 & 0.04 & 0.50 \\ 0.27 & 0.40 & 0.41 & 0.04 & 0.38 & 0.80 & 0.11 & 0.23 & 0.55\end{array}\right]$$
	$$b^{T}=[0.7,0.1,0.8,0.9,0.2,0.5,0.6]$$
	In this example, $\mathscr{I}=\{1,2, \ldots, 7\}$, $\mathscr{J}=\{1,2, \ldots, 9\}$ and $\varphi$ is the Dubois-Prade t-norm (that is continuous and non- Archimedean) with $\gamma=0.5$. So, according to Table \ref{tbl_a1} we have $\varphi(x, y)=T_{D P}^{0.5}(x, y)=x y / \max \{x, y, 0.5\}$. For $i=1$ and $j=3$, we have $a_{13}^{+}=0.8$ and $a_{13}^{-}=b_{1}=0.7$. Hence, according to Table \ref{tbl_a2}, $l_{13}^{-}=\max \left\{b_{1}, \gamma\right\}=0.7$, $l_{13}^{+}=b_{1}=0.7$, $u_{13}^{-}=1$ and $u_{13}^{+}=b_{1}=0.7$. So, from Corollary \ref{corl-1} (d) we obtain $I_{13}=\left[1-u_{13}^{-}, u_{13}^{+}\right]=[0,0.7]$ and $S_{13}=\left[1-u_{13}^{-}, 1-l_{13}^{-}\right] \cup\left[l_{13}^{+}, u_{13}^{+}\right]=[0,0.3] \cup\{0.7\}$ (since $\left.l_{i j}^{+}=0.7>0.3=1-l_{13}^{-}\right)$. Tables \ref{t-1} and \ref{t-2} show all the sets $I_{i j}$ and $S_{i j}$ for each $i \in \mathscr{I}$ and $j \in \mathscr{J}$, respectively. Moreover, by Definition \ref{def-4} and Tables \ref{t-1} and \ref{t-2}, it can be easily calculated that $I_{3}=\bigcap_{i \in \mathscr{I}} I_{i 3}=[0.1,0.7]$ and $S_{13}^{\prime}=S_{13} \cap I_{3}=[0.1,0.3] \cup\{0.7\}$. Tables \ref{t-3} and \ref{t-4} show all the sets $I_{j}$ and $S_{i j}^{\prime}$ for each $i \in \mathscr{I}$ and $j \in \mathscr{J}$, respectively.
\end{example}
% ==================================== section3_Resolution of the feasible solution set ====================================
\section{Resolution of the feasible solution set}\label{sec-3}
The following lemma gives two necessary conditions for the feasibility of Problem (\ref{eq_1}).

\begin{lemma}\label{lm-3}
	\textbf{(a)} If $S\left(A^{+}, A^{-}, b\right) \neq \varnothing$, then $I_{j} \neq \varnothing$, $\forall j \in \mathscr{J}$. \textbf{(b)} If $S\left(A^{+}, A^{-}, b\right) \neq \varnothing$, then for each $i \in \mathscr{I}$, there exists at least one $j_{i} \in \mathscr{J}$ such that $S_{ij_{i}}^{\prime} \neq \varnothing$.
\end{lemma}

\begin{proof}
	\textbf{(a)} Suppose that $x \in S\left(A^{+}, A^{-}, b\right)$ and $I_{j_{0}}=\varnothing$ for some $j_{0} \in \mathscr{J}$. Hence, there exists $i_{0} \in \mathscr{I}$ such that $x \notin I_{i_{0} j_{0}}$ (Definition \ref{def-4}). So, from Definition \ref{def-3}, we have $\max \left\{\varphi\left(a_{i_{0} j_{0}}^{+}, x_{j_{0}}\right), \varphi\left(a_{i_{0} j_{0}}^{-}, 1-x_{j_{0}}\right)\right\}>b_{i_{0}}$, that implies $x \notin S_{i_{0}}$. So, from $S\left(A^{+}, A^{-}, b\right)=\bigcap_{i \in \mathscr{I}} S_{i}$, we obtain $x \notin S\left(A^{+}, A^{-}, b\right)$ that is a contradiction. \textbf{(b)} Assume that $x \in S\left(A^{+}, A^{-}, b\right)$ and there exists some $i_{0} \in \mathscr{I}$ such that $S_{i_{0} j}^{\prime}=\varnothing$, $\forall j \in \mathscr{J}$. So, $x_{j} \notin S_{i_{0}j}^{\prime}=S_{i_{0}j} \cap I_{j}$, $\forall j \in \mathscr{J}$ (Definition \ref{def-4}). Now, if $x_{j} \notin I_{j}$ for some $j \in \mathscr{J}$, then from Part (a) we obtain $S\left(A^{+}, A^{-}, b\right) \neq \varnothing$, that is a contradiction. On the other hand, if $x_{j} \notin S_{i_{0} j}^{\prime}$ and $x_{j} \in I_{j}$, $\forall j \in \mathscr{J}$, then it is concluded that $x_{j} \notin S_{i_{0} j}$, $\forall j \in \mathscr{J}$, that implies $\max \left\{\varphi\left(a_{i_{0} j}^{+}, x_{j}\right), \varphi\left(a_{i_{0} j}^{-}, 1-x_{j}\right)\right\}<b_{i_{0}}$, $\forall j \in \mathscr{J}$ (Definition \ref{def-3}). Therefore, $x \notin S_{i_{0}}$, that contradicts the assumption that $x \in S\left(A^{+}, A^{-}, b\right)$.
\end{proof}

The following lemma provides a necessary and sufficient condition to guarantee that a given $x \in[0,1]^{n}$ is feasible to Problem (\ref{eq_1}) or not.

\begin{lemma}\label{lm-4}
	Suppose that $\varphi$ is a continuous t-norm. Then, $x \in S\left(A^{+}, A^{-}, b\right)$ if and only if the following statements hold true:
	\begin{itemize}
		\item[\textbf{(I)}] $x_{j} \in I_{j}$, $\forall j \in \mathscr{J}$.
		\item[\textbf{(II)}] For each $i \in \mathscr{I}$, there exists at least one $j_{i} \in \mathscr{J}$ such that $x_{j_{i}} \in S_{i j_{j}}^{\prime}$.
	\end{itemize}
\end{lemma}

\begin{proof}
	Suppose that $x \in[0,1]^{n}$ satisfies the conditions (I) and (II). So, from the condition (I) and Definitions \ref{def-3} and \ref{def-4}, we have $x_{j} \in I_{i j}$ and therefore $\max \left\{\varphi\left(a_{i j}^{+}, x_{j}\right), \varphi\left(a_{i j}^{-}, 1-x_{j}\right)\right\} \leq b_{i}$, $\forall i \in \mathscr{I}$ and $\forall j \in \mathscr{J}$ . On the other hand, from the condition (II), and Definitions \ref{def-3} and \ref{def-4}, for each $i \in \mathscr{I}$ we have $\max \left\{\varphi\left(a_{i j_i}^{+}, x_{j_{i}}\right), \varphi\left(a_{i j_i}^{-}, 1-x_{j_{i}}\right)\right\}=b_{i}$ for some $j_{i} \in \mathscr{J}$. Consequently, for each $i \in \mathscr{I}$, it is concluded that $\max _{j=1}^{n}\left\{\max \left\{\varphi\left(a_{i j}^{+}, x_{j}\right), \varphi\left(a_{i j}^{-}, 1-x_{j}\right)\right\}\right\}=\max \left\{\varphi\left(a_{i j_i}^{+}, x_{j_{i}}\right), \varphi\left(a_{i j_i}^{-}, 1-x_{j_{i}}\right)\right\}=b_{i}$. Thus, $x \in S_{i}$, $\forall i \in \mathscr{I}$, that implies $x \in \bigcap_{i \in \mathscr{I}} S_{i}=S\left(A^{+}, A^{-}, b\right)$. The converse statement is obtained by reversing the argument.
\end{proof}

\begin{corollary}\label{corl-2}
	For each $i \in \mathscr{I}$, $x \in S_{i}$ if and only if $x_{j} \in I_{i j}$, $\forall j \in \mathscr{J}$, and there exists at least one $j_{i} \in \mathscr{J}$ such that $x_{j_{i}} \in S_{i j_{i}}^{\prime}$.
\end{corollary}

\begin{proof}
	The proof is resulted from Lemma \ref{lm-4}, where $S\left(A^{+}, A^{-}, b\right)$ and $I_{j}$ are replaced by $S_{i}$ and $I_{i j}$, respectively.
\end{proof}

\begin{definition}\label{def-5}
	For each $i \in \mathscr{I}$, we define $\mathscr{J}_{i}=\left\{j \in \mathscr{J}: S_{i j}^{\prime} \neq \varnothing\right\}$. Similarly, for each $j \in \mathscr{J}$, define $\mathscr{I}_{j}=\left\{i \in \mathscr{J}: S_{i j}^{\prime} \neq \varnothing\right\}$.
\end{definition}

\begin{definition}\label{def-6}
	A function $e$ (on $\mathscr{I}$) is said to be an admissible function if $e(i) \in \mathscr{J}_{i}(e)$, $\forall i \in \mathscr{I}$, where
	\begin{itemize}
		\item[\textbf{(I)}] $\mathscr{J}_{1}(e)=\mathscr{J}_{1}$.
		\item[\textbf{(II)}] $\mathscr{I}_{j}(e, i)=\{k \in \mathscr{I}: 1 \leq k<i$ and $e(k)=j\}$, $\forall j \in \mathscr{J}$ and $\forall i \in \mathscr{I}-\{1\}$.
		\item[\textbf{(III)}] $\mathscr{J}_{i}(e)=\left\{j \in \mathscr{J}_{i}: \mathscr{I}_{j}(e, i)=\varnothing\right.$ or $\left.S_{i j}^{\prime} \cap\left(\bigcap_{k \in \mathscr{I}_{j}(e, i)} S_{k j}^{\prime}\right) \neq \varnothing\right\}$, $\forall i \in \mathscr{I}-\{1\}$.
	\end{itemize}
	Also, let $E$ denote the set of all the admissible functions. For the sake of convenience, we can represent each $e$ as the vector $e=\left[j_{1}, \ldots, j_{m}\right]$ in which $e(i)=j_{i}$, $\forall i \in \mathscr{I}$.
\end{definition}

\begin{definition}\label{def-7}
	For each $e \in E$, let $\mathscr{I}_{j}(e)=\{i \in \mathscr{I}: e(i)=j\}$ and $S(e)$ be the set of all the vectors $x=\left(x_{1}, \ldots, x_{n}\right)$ such that
	\begin{equation}\label{eq-5}
		x_{j} \in\left\{\begin{array}{ll}
			\bigcap_{i \in \mathscr{I}_{j}(e)} S_{i j}^{\prime} & , \mathscr{I}_{j}(e) \neq \varnothing  \\
			I_{j} & , \mathscr{I}_{j}(e)=\varnothing
		\end{array} \quad, \forall j \in \mathscr{J}\right.
	\end{equation}
	In other words, $S(e)=S(e)_{1} \times \cdots \times S(e)_{n}$ where $S(e)_{j}=\bigcap_{i \in \mathscr{I}_{j}(e)} S_{i j}^{\prime}$ if $\mathscr{I}_{j}(e) \neq \varnothing$, and $S(e)_{j}=I_{j}$ if $\mathscr{I}_{j}(e)=\varnothing$.
\end{definition}

\begin{corollary}\label{corl-3}
	Let $e: \mathscr{I} \rightarrow \bigcup_{i \in \mathscr{I}} \mathscr{J}_{i}$ be a function so that $e(i) \in \mathscr{J}_{i}$, $\forall i \in \mathscr{I}$. Then, $e \in E$ if and only if $\bigcap_{i \in \mathscr{I}_{j}(e)} S_{i j}^{\prime} \neq \varnothing$ for each $j \in \mathscr{J}$ such that $\mathscr{I}_{j}(e) \neq \varnothing$.
\end{corollary}

\begin{proof}
	Suppose that $e \in E$ and $\mathscr{I}_{j_{0}}(e) \neq \varnothing$ for some $j_{0} \in \mathscr{J}$. Also, without loss of generality, let $i_{0}=\max \mathscr{I}_{j_{0}}(e)$. So, according to Definition \ref{def-6} , we have $\mathscr{I}_{j_{0}}(e)=\mathscr{I}_{j_{0}}\left(e, i_{0}\right) \cup\left\{i_{0}\right\}$ and $S_{i_{0} j_{0}}^{\prime} \cap\left(\bigcap_{k \in \mathscr{I}_{j_{0}}\left(e, i_{0}\right)} S_{kj_{0}}^{\prime}\right) \neq \varnothing$. Thus, $\bigcap_{i \in \mathscr{I}_{j_{0}}(e)} S_{ij_{0}}^{\prime} \neq \varnothing$. To prove the converse statement, since $e(i) \in \mathscr{J}_{i}$ ($\forall i \in \mathscr{I}$), it is sufficient to show that if $e\left(i_{0}\right)=j_{0}$ and $\mathscr{I}_{j_{0}}\left(e, i_{0}\right) \neq \varnothing$, then $S_{i_{0} j_0}^{\prime} \cap\left(\bigcap_{k \in \mathscr{I}_{j_0}\left(e, i_{0}\right)} S_{kj_{0}}^{\prime}\right) \neq \varnothing$. However, since $\mathscr{I}_{j_{0}}\left(e, i_{0}\right) \subseteq \mathscr{I}_{j_{0}}(e)$, the statement $\mathscr{I}_{j_{0}}\left(e, i_{0}\right) \neq \varnothing$ implies $\mathscr{I}_{j_{0}}(e) \neq \varnothing$. Hence, from the assumption, we have $\bigcap_{i \in \mathscr{I}_{j_0}(e)} S_{ij_{0}}^{\prime} \neq \varnothing$. Therefore, since $i_{0} \in \mathscr{I}_{j_{0}}(e)$ (because $e\left(i_{0}\right)=j_{0}$) and $\mathscr{I}_{j_{0}}\left(e, i_{0}\right) \subseteq \mathscr{I}_{j_{0}}(e)$, it is concluded that $S_{i_{0} j_{0}}^{\prime} \cap\left(\bigcap_{k \in \mathscr{I}_{j_0}\left(e, i_{0}\right)} S_{kj_{0}}^{\prime}\right) \neq \varnothing$.
\end{proof}

\begin{remark}\label{rmk-2}
	According to Corollary \ref{corl-3}, the number of admissible functions is bounded above by $\prod_{i \in \mathscr{S}}\left|\mathscr{J}_{i}\right|$, where $\left|\mathscr{J}_{i}\right|$ denotes the cardinality of $\mathscr{J}_{i}$. However, the actual number of admissible functions is usually much less than this value.
\end{remark}

Through use of the admissible functions, the following theorem determines the feasible solutions set of Problem (\ref{eq_1}).

\begin{theorem}\label{thm-1}
	$S\left(A^{+}, A^{-}, b\right)=\bigcup_{e \in E} S(e)$.
\end{theorem}

\begin{proof}
	Let $x \in \bigcup_{e \in E} S(e)$. So, $x \in S\left(e_{0}\right)$ for some $e_{0} \in E$. Hence, according to (\ref{eq-5}), for each $j \in \mathscr{J}$ we have either $x_{j} \in \bigcap_{i \in \mathscr{I}_{j}(e_0)} S_{ij}^{\prime}$ (if $\mathscr{I}_{j}\left(e_{0}\right) \neq \varnothing$)
	or $x_{j} \in I_{j}$ (if $\mathscr{I}_{j}\left(e_{0}\right)=\varnothing$). But, since $S_{i j}^{\prime}=S_{i j} \cap I_{j}$ (Definition \ref{def-4}), from $x_{j} \in \bigcap_{i \in \mathscr{I}_{j}\left(e_{0}\right)} S_{ij}^{\prime}$ it is obtained again $x_{j} \in I_{j}$. Consequently, $x_{j} \in I_{j}$, $\forall j \in \mathscr{J}$. On the other hand, from Definition \ref{def-6} we have $e_{0}(i)=j_{i} \in \mathscr{J}_{i}\left(e_{0}\right) \subseteq \mathscr{J}_{i}$ ($\forall i \in \mathscr{I}$), which implies $\mathscr{I}_{j_{i}}\left(e_{0}\right) \neq \varnothing$, and therefore from (\ref{eq-5}) we have $x_{j_{i}} \in \bigcap_{k \in \mathscr{I}_{j_i}\left(e_{0}\right)} S_{kj_{i}}^{\prime} \subseteq S_{i j_{i}}^{\prime}$. Hence, $x_{j_{i}} \in S_{ij_{i}}^{\prime}$, $\forall i \in \mathscr{I}$. Now, Lemma \ref{lm-4} requires that $x \in S\left(A^{+}, A^{-}, b\right)$. Conversely, let $x \in S\left(A^{+}, A^{-}, b\right)$, $\mathscr{I}_{j}(x)=\left\{i \in \mathscr{I}: x_{j} \in S_{i j}^{\prime}\right\}$ and $\mathscr{J}_{i}(x)=\left\{j \in \mathscr{J}: x_{j} \in S_{i j}^{\prime}\right\}$. So, for each $i \in \mathscr{I}_{j}(x)$ and each $j \in \mathscr{J}_{i}(x)$, we have $S_{i j}^{\prime} \neq \varnothing$ that means $j \in \mathscr{J}_{i}$ (Definition \ref{def-5}). Also, Lemma \ref{lm-4} implies that $\mathscr{J}_{i}(x) \neq \varnothing$, $\forall i \in \mathscr{I}$. Without loss of generality, let $j_{i}=\min \mathscr{J}_{i}(x)$ and $e_{0}(i)=j_{i}$, $\forall i \in \mathscr{I}$. Therefore, $e_{0}$ is a function on $\mathscr{I}$ such that
	\begin{equation}\label{eq-6}
		e_{0}(i)=j_{i} \in \mathscr{J}_{i}, ~\forall i \in \mathscr{I} 
	\end{equation}
	Moreover, we have
	\begin{equation}\label{eq-7}
		\mathscr{I}_{j}\left(e_{0}\right) \neq \varnothing, ~\forall j \in\left\{j_{1}, \ldots, j_{m}\right\} ~\text { and } ~ \mathscr{I}_{j}\left(e_{0}\right)=\varnothing, ~\forall j \in \mathscr{J}-\left\{j_{1}, \ldots, j_{m}\right\} 
	\end{equation}
	In addition, if $j_{p} \in\left\{j_{1}, \ldots, j_{m}\right\}$ and $k \in \mathscr{I}_{j_{p}}\left(e_{0}\right)$ (i.e., $e_{0}(k)=j_{p}$), then by our definition it follows that $j_{p}=\min \mathscr{J}_{k}(x)$ which means $x_{j_{p}} \in S_{k j_{p}}^{\prime}$. Hence, it is concluded that
	\begin{equation}\label{eq-8}
		x_{j} \in \bigcap_{i \in \mathscr{I}_{j}\left(e_{0}\right)} S_{i j}^{\prime}, ~\forall j \in\left\{j_{1}, \ldots, j_{m}\right\}
	\end{equation}
	Consequently, by Corollary \ref{corl-3} and (\ref{eq-6}) - (\ref{eq-8}) we have $e_{0} \in E$. Also, since $x \in S\left(A^{+}, A^{-}, b\right)$, then $x_{j} \in I_{j}$, $\forall j \in \mathscr{J}$ (Lemma \ref{lm-4}). Particularly, $x_{j} \in I_{j}$, if $\mathscr{I}_{j}\left(e_{0}\right)=\varnothing$. This fact together with (\ref{eq-5}), (\ref{eq-7}) and (\ref{eq-8}) imply $x \in S\left(e_{0}\right)$.
\end{proof}

Suppose that $e \in E$, $x \in S(e)$ and $j \in \mathscr{J}$. So, According to (\ref{eq-5}), we have either $x_{j} \in \bigcap_{i \in \mathscr{I}_{j}(e)} S_{i j}^{\prime}$ (if $\mathscr{I}_{j}(e) \neq \varnothing$) or $x_{j} \in I_{j}=\left[L_{j}, U_{j}\right]$ (if $\mathscr{I}_{j}(e)=\varnothing$). On the other hand, in the former case, we have $\bigcap_{i \in \mathscr{I}_{j}(e)} S_{i j}^{\prime} \neq \varnothing$ (from Corollary \ref{corl-3}). Consequently, by noting Remark 1, it follows that $\bigcap_{i \in \mathscr{I}_{j}(e)} S_{i j}^{\prime}$ is a closed and bounded (but, not necessarily connected). So, it is concluded that $S(e)$ ($\forall e \in E$) is indeed a closed bounded subset (and therefore, a compact subset) of $[0,1]^{n}$. This fact together with Theorem \ref{t-1} imply that the feasible solutions set of Problem (\ref{eq_1}) is determined by the union of a finite number of compact sets, which are not necessarily connected.

\begin{example}\label{ex-2}
	Consider the problem stated in Example \ref{ex-1}. According to Definition \ref{def-5} and Table \ref{t-4}, we obtain $\mathscr{J}_{1}=\{3,5\}$, $\mathscr{J}_{2}=\{1,2,7\}$, $\mathscr{J}_{3}=\{8,9\}$, $\mathscr{J}_{4}=\{2,3,4,7\}$, $\mathscr{J}_{5}=\{5\}$, $\mathscr{J}_{6}=\{8,9\}$ and $\mathscr{J}_{7}=\{2,6\}$. Hence, according to Remark \ref{rmk-2}, the number of admissible functions is bounded above by $\prod_{i \in \mathscr{I}}\left|\mathscr{J}_{i}\right|=2 \times 3 \times 2 \times 4 \times 1 \times 2 \times 2=192$. Now, noting Definition \ref{def-6} and Corollary \ref{corl-3}, consider functions $e_{1}=[3,2,9,2,5,9,6]$ and $e_{2}=[3,2,9,3,5,9,2]$ from $\mathscr{I}$ to $\bigcup_{i \in \mathscr{I}} \mathscr{J}_{i}$ so that $e_{p}(i) \in \mathscr{J}_{i}$ for any $i \in \mathscr{I}$ and $p \in\{1,2\}$. So, $e_{1}(2)=e_{1}(4)=2$ and $\mathscr{I}_{2}\left(e_{1}, 4\right)=\{2\}$ (Definition \ref{def-6} (II)). Also, $S_{22}^{\prime}=\{0.75\}$ and $S_{42}^{\prime}=\{0.9\}$ (see Table \ref{t-4}). Therefore, $\mathscr{I}_{2}\left(e_{1}, 4\right) \neq \varnothing$ and $S_{42}^{\prime} \cap\left(\bigcap_{k \in \mathscr{I}_{2}\left(e_{1}, 4\right)} S_{k 2}^{\prime}\right)=S_{42}^{\prime} \cap S_{22}^{\prime}=\varnothing$, which implies $2 \notin \mathscr{J}_{4}\left(e_{1}\right)$ (Definition \ref{def-6} ( III)). As a result, since $e_{1}(4)=2 \notin \mathscr{J}_{4}\left(e_{1}\right)$, based on Definition \ref{def-6} it follows that $e_{1}$ is not admissible. However, $e_{2}$ is indeed an admissible function. More precisely, we have $\mathscr{I}_{2}\left(e_{2}, 2\right)=\mathscr{I}_{9}\left(e_{2}, 3\right)=\mathscr{I}_{4}\left(e_{2}, 4\right)=\mathscr{I}_{5}\left(e_{2}, 5\right)=\mathscr{I}_{8}\left(e_{2}, 6\right)=\varnothing$ and $\mathscr{I}_{2}\left(e_{2}, 7\right)=\{2\}$. Also, for the non-empty set $\mathscr{I}_{2}\left(e_{2}, 7\right)$, we have $S_{72}^{\prime} \cap\left(\bigcap_{k \in \mathscr{I}_{2}(e, 7)} S_{k 2}^{\prime}\right)=S_{72}^{\prime} \cap S_{22}^{\prime}=\{0.75\} \neq \varnothing$. So, the three conditions of Definition \ref{def-6} hold for $e_{2}$. Moreover, by noting Table \ref{t-4}, it follows that $\bigcap_{i \in \mathscr{I}_{1}\left(e_{2}\right)} S_{i 1}^{\prime}=\bigcap_{i \in \mathscr{I}_{6}\left(e_{2}\right)} S_{i6}^{\prime}=\bigcap_{i \in \mathscr{I}_{7}\left(e_{2}\right)} S_{i 7}^{\prime}=\varnothing$, $\bigcap_{i \in \mathscr{I}_{2}\left(e_{2}\right)} S_{i 2}^{\prime}=S_{22}^{\prime} \cap S_{72}^{\prime}=\{0.75\}$, $\bigcap_{i \in \mathscr{I}_{3}\left(e_{2}\right)} S_{i 3}^{\prime}=S_{13}^{\prime}=[0.1,0.3] \cup\{0.7\}$, $\bigcap_{i \in \mathscr{I}_{4}\left(e_{2}\right)} S_{i 4}^{\prime}=S_{44}^{\prime}=[0,0.1] \cup[0.9,1]$, $\bigcap_{i \in \mathscr{I}_{5}\left(e_{2}\right)} S_{i5}^{\prime}=S_{55}^{\prime}=\{0.75\}$, $\bigcap_{i \in \mathscr{I}_{8}\left(e_{2}\right)} S_{i8}^{\prime}=S_{68}^{\prime}=[0.5,1]$ and $\bigcap_{i \in \mathscr{I}_{9}\left(e_{2}\right)} S_{i9}^{\prime}=S_{39}^{\prime}=\{0.2\}$. Furthermore, from Table 3, $I_{1}=[0,0.25]$, $I_{6}=[0.4,0.6]$ and $I_{7}=\{0.1\}$. Hence, from Definition \ref{def-7},$ S\left(e_{2}\right)$ is calculated as the Cartesian product
	$S\left(e_{2}\right)=[0,0.25] \times\{0.75\} \times[0.1,0.3] \cup\{0.7\} \times[0,0.1] \cup[0.9,1] \times\{0.75\} \times[0.4,0.6] \times\{0.1\} \times[0.5,1] \times\{0.2\}$.
\end{example}

% ==================================== section4_Simplification rules ====================================
\section{Simplification rules}\label{sec-4}
This subsection describes some simplification techniques that can be used to accelerate the determination of the feasible region. Throughout this section, $S_{\{i\}}\left(A^{+}, A^{-}, b\right)$ denotes the feasible region of the reduced problem obtained by removing the $i$'th equation from Problem (\ref{eq_1}). The proofs of the following rules have been presented in Appendix $B$.

\begin{therule}\label{rule-1}
	Suppose that $S\left(A^{+}, A^{-}, b\right) \neq \varnothing$ and $i_{0} \in \mathscr{I}$. If $b_{i_0}=0$, then the $i_{0}$'th equation is a redundant constraint and it can be deleted.
\end{therule}

\begin{therule}\label{rule-2}
	Suppose that $S\left(A^{+}, A^{-}, b\right) \neq \varnothing$ and there exists $j_{0} \in \mathscr{J}$ such that $I_{j_{0}}=\{k\}$ is a singleton set. Then, $x_{j_{0}}=k$ for each feasible solution $x$. Also, the $j_{0}$'th column and any equation $i_{0}$ such that $k \in S_{i_0 j_0}^{\prime}$ can be removed from the problem.
\end{therule}

\begin{therule}\label{rule-3}
	Suppose that $S\left(A^{+}, A^{-}, b\right) \neq \varnothing$ and there exist $i, i_{0} \in \mathscr{I}$ such that $S_{i j}^{\prime} \subseteq S_{i_0 j}^{\prime}$, $\forall j \in \mathscr{J}$. Then, the $i_{0}$'th equation is a redundant constraint and it can be deleted.
\end{therule}

\begin{therule}\label{rule-4}
	Suppose that $S\left(A^{+}, A^{-}, b\right) \neq \varnothing$ and there exist $i_{0} \in \mathscr{I}$ and $j_{0} \in \mathscr{J}$ such that $\mathscr{J}_{i_0}=\left\{j_{0}\right\}$ and $S_{i_0 j_0}^{\prime}=\{k\}$ are singleton sets. Then, $x_{j_{0}}=k$ for each feasible solution $x$. Also, the $j_{0}$'th column and any equation $i$ such that $k \in S_{ij_{0}}^{\prime}$ can be removed from the problem.
\end{therule}

\begin{therule}\label{rule-5}
	Suppose that $S\left(A^{+}, A^{-}, b\right) \neq \varnothing$ and there exist $i_{0} \in \mathscr{I}$ and $j_{0} \in \mathscr{J}$ such that $S_{i_0 j_0}^{\prime}=I_{j_{0}}$. Then, the $i_{0}$'th equation is a redundant constraint and it can be deleted.
\end{therule}

We now summarize the preceding discussion as an algorithm.

\begin{algorithm}[!ht]\label{myalgorithm}
	\DontPrintSemicolon
	\SetKwInOut{Input}{Input}
	\SetKwInOut{Output}{Output}
	
	\Input{Given Problem (\ref{eq_1}):}
	\Output{Generate $S$ (the feasible solution set)}
	
	Compute sets $I_{i j}$, $S_{i j}$, $I_{j}$ and $S_{i j}^{\prime}$ for each $i \in \mathscr{I}$ and each $j \in \mathscr{J}$ (Corollary \ref{corl-1}, Definition \ref{def-4} and Table \ref{tbl_a2}).
	
	\If{$I_{j}=\varnothing$ for some $j \in \mathscr{J}$,}{
		stop; the problem is infeasible (Lemma \ref{lm-3} (a)).
	}
	\Else{
		\If{$S_{i j}^{\prime}=\varnothing$ for some $i \in \mathscr{I}$ and each $j \in \mathscr{J}$,}{
			the problem is infeasible (Lemma \ref{lm-3} (b)).
		}
		\Else{
			Apply the simplification rules \ref{rule-1} - \ref{rule-5} to determine the values of as many variables as possible. Delete the corresponding columns and redundant equations from the problem.
			
			Find the set $S(e)$ for each admissible function $e \in E$ in the remaining problem (Definitions \ref{def-6} and \ref{def-7}). Add the assigned variables (obtained in Step \ref{st-4}) to $S(e)$.
			
			Generate the feasible solution set $S\left(A^{+}, A^{-}, b\right)$ by $\bigcup_{e \in E} S(e)$ (Theorem \ref{thm-1}).
		}
	}
	\caption{Resolution of Problem (\ref{eq_1})}
\end{algorithm}

It is to be noted that the FRE problems whose feasible regions are defined by the constraints $A \varphi x = b$ (i.e., $\max _{j=1}^{n}\{\varphi(a_{i j}$, $x_{j})\}=b_{i}$, $\forall i \in \mathscr{I}$), are indeed special cases of Problem (\ref{eq_1}). To see this, let $A^{+}=A$ and $A^{-}=\mathbf{0}_{m \times n}$ where $\mathbf{0}_{m \times n}$ is the zero matrix. So, $a_{i j}^{+}=a_{i j}$ and $a_{i j}^{-}=0$ ($\forall i \in \mathscr{I}$ and $\forall j \in \mathscr{J}$), and therefore each constraint $\max _{j=1}^{n}\left\{\max \left\{\varphi\left(a_{i j}^{+}, x_{j}\right), \varphi\left(a_{i j}^{-}, 1-x_{j}\right)\right\}\right\}=b_{i}$ is reduced to $\max _{j=1}^{n}\left\{\max \left\{\varphi\left(a_{i j}, x_{j}\right), \varphi\left(0,1-x_{j}\right)\right\}\right\}=b_{i}$. Furthermore, by the identity law of $\varphi$, it follows that $\varphi\left(0,1-x_{j}\right)=0$, $\forall j \in J$. Consequently, the above constraints are converted into $\max _{j=1}^{n}\left\{\varphi\left(a_{i j}, x_{j}\right)\right\}=b_{i}$, $\forall i \in \mathscr{I}$. As a result, Algorithm 1 can be also used to solve the FRE problems. However, although the above algorithm accelerates the resolution of the problem by taking advantage of the five simplification rules, it should be noted that the number of $S(e)$ ($e \in E$) may grow exponentially in terms of the problem size. Consequently, finding other simplification rules or providing some methods (such as appropriate branch-and-bound techniques) is always in demand.

\begin{example}\label{ex-3}
	Let us look at the problem stated in Examples \ref{ex-1} and \ref{ex-2}. The steps of Algorithm 1 are as follows:
	
	\paragraph{Step 1}\label{st-1}
	{Tables \ref{t-1} and \ref{t-2} show all the sets $I_{i j}$ and $S_{i j}$ for each $i \in \mathscr{I}$ and $j \in \mathscr{J}$, respectively. In both tables, row $i$ ($i \in \mathscr{I}$) corresponds to equation $i$ and column $j$ ($j \in \mathscr{J}$) corresponds to variable $x_{j}$. By Table \ref{t-1} and Definition \ref{def-4}, all the sets $I_{j}$ are shown in Table \ref{t-3} for each $j \in \mathscr{J}$. Also, all the sets $S_{i j}^{\prime}$ ($\forall i \in \mathscr{I}$ and $\forall j \in \mathscr{J}$) are summarized in Table \ref{t-4} by using Table \ref{t-2}, Table \ref{t-3} and Definition \ref{def-4}.
		\begin{table}[!ht]
			\centering
			\caption{Sets $I_{i j}$ for each $i \in \mathscr{I}$ and $j \in \mathscr{J}$.}
			\label{t-1}
			\scalebox{0.99}{
				\begin{tabular}{|c|c|c|c|c|c|c|c|c|}
					\hline
					$[0,1]$ & $[0,1]$ & $[0,0.7]$ & $[0,1]$ & $[0.3,1]$ & $[0,1]$ & $[0,1]$ & $[0,1]$ & $[0,1]$ \\
					\hline
					$[0,0.25]$ & $[0.75,1]$ & $[0,1]$ & $[0,1]$ & $[0,1]$ & $[0,1]$ & $[0,0.1]$ & $[0,1]$ & $[0,1]$ \\
					\hline
					$[0,1]$ & $[0,1]$ & $[0,1]$ & $[0,1]$ & $[0,1]$ & $[0,1]$ & $[0,1]$ & $[0,1]$ & $[0.2,1]$ \\
					\hline
					$[0,1]$ & $[0,0.9]$ & $[0.1,1]$ & $[0,1]$ & $[0,1]$ & $[0,1]$ & $[0.1,1]$ & $[0,1]$ & $[0,1]$ \\
					\hline
					$[0,1]$ & $[0.6,1]$ & $[0,1]$ & $[0,1]$ & $[0.75,1]$ & $[0,1]$ & $[0,1]$ & $[0,1]$ & $[0,1]$ \\
					\hline
					$[0,1]$ & $[0,1]$ & $[0,1]$ & $[0,1]$ & $[0,1]$ & $[0,1]$ & $[0,1]$ & $[0,1]$ & $[0,1]$ \\
					\hline
					$[0,1]$ & $[0,1]$ & $[0,1]$ & $[0,1]$ & $[0,1]$ & $[0.4,0.6]$ & $[0,1]$ & $[0,1]$ & $[0,1]$ \\
					\hline
			\end{tabular}}
		\end{table}
		
		\begin{table}[!ht]
			\centering
			\caption{Sets $S_{i j}$ for each $i \in \mathscr{I}$ and $j \in \mathscr{J}$.}
			\label{t-2}
			\scalebox{0.85}{
				\begin{tabular}{|c|c|c|c|c|c|c|c|c|}
					\hline
					$\varnothing$ & $\varnothing$ & \begin{tabular}{c}
						$[0,0.3] \cup \{0.7\}$ \\
					\end{tabular} & $\varnothing$ & \begin{tabular}{c}
						$\{0.3\} \cup [0.7,1]$ \\
					\end{tabular} & $\varnothing$ & $\varnothing$ & $\varnothing$ & $\varnothing$ \\
					\hline
					$[0,0.25]$ & $\{0.75\}$ & $\varnothing$ & $\varnothing$ & $\varnothing$ & $\varnothing$ & $\{0.1\}$ & $\varnothing$ & $\varnothing$ \\
					\hline
					$\varnothing$ & $\varnothing$ & $\varnothing$ & $\varnothing$ & $\varnothing$ & $\varnothing$ & $\varnothing$ & \begin{tabular}{l}
						$[0,0.2] \cup [0.8,1]$ \\
					\end{tabular} & $\{0.2\}$ \\
					\hline
					$\varnothing$ & $\{0.9\}$ & $\{0.1\}$ & \begin{tabular}{l}
						$[0,0.1] \cup [0.9,1]$ \\
					\end{tabular} & $\varnothing$ & $\varnothing$ & $\{0.1\}$ & $\varnothing$ & $\varnothing$ \\
					\hline
					$\varnothing$ & $\{0.6\}$ & $\varnothing$ & $\varnothing$ & $\{0.75\}$ & $\varnothing$ & $\varnothing$ & $\varnothing$ & $\varnothing$ \\
					\hline
					$\varnothing$ & $\varnothing$ & $\varnothing$ & $\varnothing$ & $\varnothing$ & $\varnothing$ & $\varnothing$ & $[0.5,1]$ & $[0,0.5]$ \\
					\hline
					$\varnothing$ & $[0.6,1]$ & $\varnothing$ & $\varnothing$ & $\varnothing$ & \begin{tabular}{l} $\{0.4\} \cup \{0.6\}$ \\ \end{tabular} & $\varnothing$ & $\varnothing$ & $\varnothing$ \\
					\hline
			\end{tabular}}
		\end{table}
		
		\begin{table}[!ht]
			\centering
			\caption{Sets $I_{j}=\left[L_{j}, U_{j}\right]$ for each $j \in \mathscr{J}$.}
			\label{t-3}
			\scalebox{0.99}{
				\begin{tabular}{|l|l|l|l|l|l|l|l|l|}
					\hline
					$[0,0.25]$ & $[0.75,0.9]$ & $[0.1,0.7]$ & $[0,1]$ & $[0.75,1]$ & $[0.4,0.6]$ & $\{0.1\}$ & $[0,1]$ & $[0.2,1]$ \\
					\hline
			\end{tabular}}
		\end{table}
		
		\begin{table}[!ht]
			\centering
			\caption{Sets $S_{i j}^{\prime}$ for each $i \in \mathscr{I}$ and $j \in \mathscr{J}$.}
			\label{t-4}
			\scalebox{0.85}{
				\begin{tabular}{|c|c|c|c|c|c|c|c|c|}
					\hline
					$\varnothing$ & $\varnothing$ & \begin{tabular}{c}
						$[0.1,0.3] \cup \{0.7\}$ \\
					\end{tabular} & $\varnothing$ & $[0.75,1]$ & $\varnothing$ & $\varnothing$ & $\varnothing$ & $\varnothing$ \\
					\hline
					$[0,0.25]$ & $\{0.75\}$ & $\varnothing$ & $\varnothing$ & $\varnothing$ & $\varnothing$ & $\{0.1\}$ & $\varnothing$ & $\varnothing$ \\
					\hline
					$\varnothing$ & $\varnothing$ & $\varnothing$ & $\varnothing$ & $\varnothing$ & $\varnothing$ & $\varnothing$ & \begin{tabular}{c}
						$[0,0.2] \cup [0.8,1]$ \\
					\end{tabular} & $\{0.2\}$ \\
					\hline
					$\varnothing$ & $\{0.9\}$ & $\{0.1\}$ & \begin{tabular}{l}
						$[0,0.1] \cup [0.9,1]$ \\
					\end{tabular} & $\varnothing$ & $\varnothing$ & $\{0.1\}$ & $\varnothing$ & $\varnothing$ \\
					\hline
					$\varnothing$ & $\varnothing$ & $\varnothing$ & $\varnothing$ & $\{0.75\}$ & $\varnothing$ & $\varnothing$ & $\varnothing$ & $\varnothing$ \\
					\hline
					$\varnothing$ & $\varnothing$ & $\varnothing$ & $\varnothing$ & $\varnothing$ & $\varnothing$ & $\varnothing$ & $[0.5,1]$ & $[0.2,0.5]$ \\
					\hline
					$\varnothing$ & $[0.75,0.9]$ & $\varnothing$ & $\varnothing$ & $\varnothing$ & \begin{tabular}{l}
						$\{0.4\} \cup\{0.6\}$ \\
					\end{tabular} & $\varnothing$ & $\varnothing$ & $\varnothing$ \\
					\hline
			\end{tabular}}
		\end{table}
	}
	
	\paragraph{Step 2}\label{st-2}
	{From Table \ref{t-3}, it is clear that $I_{j} \neq \varnothing$, $\forall j \in \mathscr{J}$.}
	
	\paragraph{Step 3}\label{st-3}
	{As it was described in Table \ref{t-4}, for each $i \in \mathscr{I}$ there exists at least one $j_{i} \in \mathscr{J}$ such that $S_{{ij}_{i}}^{\prime} \neq \varnothing$.}
	
	\paragraph{Step 4}\label{st-4}
	{Since $I_{7}=\{0.1\}$ is a singleton set (see Table \ref{t-3}) and $0.1 \in S_{27}^{\prime} \cap S_{47}^{\prime}$ (see Table \ref{t-4}), Rule \ref{rule-2} indicates that $x_{7}=0.1$ for each feasible solution $x$ (particularly, $x_{7}^{*}=0.1$ for each optimal solution $x^{*}$). Also, column 7 and both rows 2 and 4 can be removed from the problem. So, by this simplification rule, the upper bound of the number of admissible functions is reduced from 192 (Example \ref{ex-2}) to $\prod_{i \in \mathscr{I}-(2,4\}}\left|\mathscr{J}_{i}\right|=2 \times 2 \times 1 \times 2 \times 2=16$; that is $|E| \leq 16$. Subsequently, according to Table \ref{t-4}, it turns out that $\mathscr{J}_{5}=\{5\}$ and $S_{55}^{\prime}=\{0.75\}$ are singleton sets and $0.75 \in S_{15}^{\prime}$. Thus, $x_{5}$ is assigned to $x_{5}=0.75$, and also column 5 and rows 1 and 5 can be deleted by Rule \ref{rule-4}. Consequently, we have $\prod_{i \in \mathscr{I}-(1,2,4,5)}\left|\mathscr{J}_{i}\right|=2 \times 2 \times 2=8$ and $|E| \leq 8$. It is to be noted that the first row can be also deleted by Rule \ref{rule-3}; because, $S_{5 j}^{\prime} \subseteq S_{1 j}^{\prime}$, $\forall j \in \mathscr{J}$. Moreover, as shown in Tables \ref{t-3} and \ref{t-4}, $S_{72}^{\prime}=I_{2}$. So, it is concluded that the seventh equation is also a redundant constraint and it can be deleted by Rule \ref{rule-5}. Hence, by applying this simplification technique, the upper bound of the number of admissible functions is further reduced from $8$ to $|E|=\prod_{i \in \mathscr{I}-\{1,2,4,5,7\}}\left|\mathscr{J}_{i}\right|=2 \times 2=4$. So, after applying the above simplification rules, columns $\{5,7\}$ and rows $\{1,2,4,5,7\}$ are deleted and we obtain $x_{5}^{*}=0.75$ and $x_{7}^{*}=0.1$. Therefore, the reduced matrices $A^{+}$and $A^{-}$, and the right-hand-side vector $b$ become
		$$
		\begin{gathered}
			A^{+}=\left[\begin{array}{lllllll}
				0.72 & 0.23 & 0.75 & 0.44 & 0.61 & 0.80 & 0.67 \\
				0.28 & 0.43 & 0.35 & 0.28 & 0.22 & 0.50 & 0.00
			\end{array}\right] ,\quad A^{-}=\left[\begin{array}{lllllll}
				0.13 & 0.63 & 0.74 & 0.25 & 0.73 & 0.80 & 0.90 \\
				0.00 & 0.29 & 0.33 & 0.47 & 0.34 & 0.04 & 0.50
			\end{array}\right] \\
			b^{T}=[0.8,0.5]
		\end{gathered}
		$$
		The current matrices $A^{+}$and $A^{-}$are equivalent to two rows (rows 3 and 6 in the main problem) and seven columns (columns $1-4,6,8$ and 9 in the main problem). Furthermore, Table \ref{t-4} is updated as follows:
		
		\begin{table}[!ht]
			\centering
			\caption{Updated Table \ref{t-4} after applying Rules \ref{rule-1} - \ref{rule-5}.}
			\label{t-5}
			\scalebox{0.99}{
				\begin{tabular}{|c|c|c|c|c|c|c|}
					\hline
					$\varnothing$ & $\varnothing$ & $\varnothing$ & $\varnothing$ & $\varnothing$ & $[0,0.2] \cup[0.8,1]$ & $\{0.2\}$ \\
					\hline
					$\varnothing$ & $\varnothing$ & $\varnothing$ & $\varnothing$ & $\varnothing$ & $[0.5,1]$ & $[0.2,0.5]$ \\
					\hline
			\end{tabular}}
		\end{table}
	}
	
	\paragraph{Step 5 and 6}\label{st-56}
	{Based on the results obtained in the previous step, we have $|E|=4$. Four admissible functions $e \in E$ are obtained from Table \ref{t-5} and Definition \ref{def-6} as $e_{1}=[6,6]$, $e_{2}=[6,7]$, $e_{3}=[7,6]$ and $e_{4}=[7,7]$. So, by Definition \ref{def-7}, their corresponding sets $S(e)$ are as follows:
		$S\left(e_{1}\right)=[0,1] \times[0,1] \times[0,1] \times[0,1] \times[0,1] \times[0.8,1] \times[0,1]$
		
		$S\left(e_{2}\right)=[0,1] \times[0,1] \times[0,1] \times[0,1] \times[0,1] \times[0,0.2] \cup[0.8,1] \times[0.2,0.5]$
		
		$S\left(e_{3}\right)=[0,1] \times[0,1] \times[0,1] \times[0,1] \times[0,1] \times[0.5,1] \times\{0.2\}$
		
		$S\left(e_{4}\right)=[0,1] \times[0,1] \times[0,1] \times[0,1] \times[0,1] \times[0,1] \times\{0.2\}$
		
		Now, by adding the assigned variables $x_{5}=0.75$ and $x_{7}=0.1$, we obtain
		
		$S\left(e_{1}\right)=[0,1] \times[0,1] \times[0,1] \times[0,1] \times\{0.75\} \times[0,1] \times\{0.1\} \times[0.8,1] \times[0,1] $
		
		$S\left(e_{2}\right)=[0,1] \times[0,1] \times[0,1] \times[0,1] \times\{0.75\} \times[0,1] \times\{0.1\} \times[0,0.2] \cup[0.8,1] \times[0.2,0.5]$
		
		$S\left(e_{3}\right)=[0,1] \times[0,1] \times[0,1] \times[0,1] \times\{0.75\} \times[0,1] \times\{0.1\} \times[0.5,1] \times\{0.2\}$
		
		$S\left(e_{4}\right)=[0,1] \times[0,1] \times[0,1] \times[0,1] \times\{0.75\} \times[0,1] \times\{0.1\} \times[0,1] \times\{0.2\}$
		
		Therefore, from Theorem $1, S\left(A^{+}, A^{-}, b\right)=\bigcup_{i=1}^{4} S\left(e_{i}\right)$.}
\end{example}

% ==================================== section5_Local and global optimal solutions ====================================
\section{Local and global optimal solutions}\label{sec-5}
As mentioned before, throughout this section we assume $f$ is a continuous function, that is non-decreasing (non-increasing) in $x_{j}$, $\forall j \in \mathscr{J}^{+}$ ($\forall j \in \mathscr{J}^{-}$).

\begin{definition}\label{def-8}
	Suppose that $S\left(A^{+}, A^{-}, b\right) \neq \varnothing$. For each $e \in E$, we define $x^{*}(e)=\left(x^{*}(e)_{1}, \ldots, x^{*}(e)_{n}\right)$ where for each $j \in \mathscr{J}$, components $x^{*}(e)_{j}$ are defined as follows:
	\begin{equation}\label{eq-9}
		x^{*}(e)_{j}= \begin{cases}\min \left\{\bigcap_{i \in \mathscr{I}_{j}(e)} S_{i j}^{\prime}\right\} & , j \in \mathscr{J}^{+} ~\text {and } \mathscr{I}_{j}(e) \neq \varnothing  \\ L_{j} & , j \in \mathscr{J}^{+} ~\text {and } \mathscr{I}_{j}(e)=\varnothing \\ \max \left\{\bigcap_{i \in \mathscr{I}_{j}(e)} S_{i j}^{\prime}\right\} & , j \in \mathscr{J}^{-} ~\text {and } \mathscr{I}_{j}(e) \neq \varnothing \\ U_{j} & , j \in \mathscr{J}^{-} ~\text {and } \mathscr{I}_{j}(e)=\varnothing\end{cases}
	\end{equation}
	where $L_{j}$ and $U_{j}$ are the lower and upper bounds of $I_{j}=\left[L_{j}, U_{j}\right]$. Also, define $F^{*}=\left\{x^{*}(e): e \in E\right\}$.
\end{definition}

According to the following theorem, any $x^{*}(e)$ is a feasible local optimal solution of problem (\ref{eq_1}).

\begin{theorem}\label{thm-2}
	Suppose that $S\left(A^{+}, A^{-}, b\right) \neq \varnothing$.
	\begin{itemize}
		\item[\textbf{(a)}] $F^{*} \subseteq S\left(A^{+}, A^{-}, b\right)$.
		\item[\textbf{(b)}] $x^{*}(e)$ is a global optimum in $S(e)$
	\end{itemize}
\end{theorem}

\begin{proof}
	\textbf{(a)} Let $x^{*}\left(e_{0}\right) \in F^{*}$ (for some $e_{0} \in E$). So, from (\ref{eq-5}) and (\ref{eq-9}), it follows that $x^{*}\left(e_{0}\right) \in S\left(e_{0}\right) \subseteq \bigcup_{e \in E} S(e)$. Hence, by Theorem 1, it is concluded that $F^{*} \subseteq S\left(A^{+}, A^{-}, b\right)$. \textbf{(b)} Let $x \in S(e)$. So, from (5) and (9), we have $x^{*}(e)_{j} \leq x_{j}$, $\forall j \in \mathscr{J}^{+}$, and $x_{j} \leq x^{*}(e)_{j}$, $\forall j \in \mathscr{J}^{-}$. Hence, $f\left(x^{*}(e)\right) \leq f(x)$ for any $x \in S(e)$.
\end{proof}

\begin{theorem}\label{thm-3}
	Suppose that $S\left(A^{+}, A^{-}, b\right) \neq \varnothing$ and $S^{*}$ denotes the set of optimal solutions for Problem (\ref{eq_1}). If $f\left(x^{*}\left(e^{*}\right)\right)=\min \left\{x^{*}(e): e \in E\right\}$, then $x^{*}\left(e^{*}\right)$ is a global optimal solution of Problem (\ref{eq_1}).
\end{theorem}

\begin{proof}
	Let $x \in S\left(A^{+}, A^{-}, b\right)$ be an arbitrary feasible solution. We shall show that $f\left(x^{*}\left(e^{*}\right)\right) \leq f(x)$. From Theorem \ref{thm-1}, $x \in S\left(e_{0}\right)$ for some $e_{0} \in E$. So, theorem \ref{thm-2} requires $f\left(x^{*}\left(e_{0}\right)\right) \leq f(x)$, which together with $f\left(x^{*}\left(e^{*}\right)\right) \leq f\left(x^{*}\left(e_{0}\right)\right)$ imply $f\left(x^{*}\left(e^{*}\right)\right) \leq f(x)$.
\end{proof}

From Theorem \ref{thm-2}, $x^{*}(e)$ is a local optimal solution, $\forall e \in E$. So, Theorem \ref{thm-3} provides also a necessary optimality condition in the sense that if $x^{*}$ is an optimal solution of Problem (\ref{eq_1}), then it must belong to $F^{*}$. In other words, each solution $x^{*}(e) \in F^{*}$ is an optimal candidate solution. Additionally, the global optimal value is obtained as $f\left(x^{*}\left(e^{*}\right)\right)=\min \left\{x^{*}(e): e \in E\right\}$.

\begin{example}[\textbf{Linear functions}]\label{ex-4}
	Assume that the objective function of Problem (\ref{eq_1}) is defined as the linear form $f(x)=\sum_{j=1}^{n} c_{j} x_{j}$, where $c_{j} \in \mathbb{R}$, $\forall j \in \mathscr{J}$. So, it is clear that $\mathscr{J}^{+}=\left\{j \in \mathscr{J}: c_{j} \geq 0\right\}$ and $\mathscr{J}^{-}=\left\{j \in \mathscr{J}: c_{j}<0\right\}$. For instance, consider the problem stated in Examples \ref{ex-1},\ref{ex-2} and \ref{ex-3} with the following objective function:
	$$
	f_{1}(x)=2 x_{1}+x_{2}-x_{3}-5 x_{4}+x_{5}+3 x_{6}-x_{7}+4 x_{8}-x_{9}
	$$
	In this example, we have $\mathscr{J}^{+}=\{1,2,5,6,8\}$ and $\mathscr{J}^{-}=\{3,4,7,9\}$. To calculate $x^{*}\left(e_{1}\right)$, we note that $e_{1}=[6,6]$, and therefore $e_{1}(1)=6$ and $e_{1}(2)=6$ (see Table \ref{t-5}). However, by noting that the sixth column in Table \ref{t-5} corresponds to the eighth column of the main problem, it follows that $\mathscr{I}_{8}\left(e_{1}\right)=\{1,2\}$. Moreover, $\mathscr{I}_{j}\left(e_{1}\right)=\varnothing$, $\forall j \in\{1, \ldots, 7\} \bigcup\{9\}$. Hence, from Definition \ref{def-8} (relation (\ref{eq-9})), it is concluded that
	$$
	x^{*}\left(e_{1}\right)_{8}  =\min \left\{\bigcap_{i \in \mathscr{I}_{8}(e)} S_{i 8}^{\prime}\right\}=\min \left\{S_{18}^{\prime} \cap S_{28}^{\prime}\right\} 
	=\min \{\{[0,0.2] \cup[0.8,1]\} \cap[0.5,1]\}=\min \{[0.8,1]\}=0.8
	$$
	$x^{*}\left(e_{1}\right)_{j}=L_{j}$, $\forall j \in\{1,2,5,6\}$, and $x^{*}\left(e_{1}\right)_{j}=U_{j}$, $\forall j \in\{3,4,7,9\}$ (see Table \ref{t-3}); that is, $$x^{*}\left(e_{1}\right)=[0,0.75,0.7,1,0.75,0.4,0.1,0.8,1].$$ By the similar calculation, we obtain:
	
	$x^{*}(e_{2})=[0,0.75,0.7,1,0.75,0.4,0.1,0,0.5]$
	
	$x^{*}(e_{3})=[0,0.75,0.7,1,0.75,0.4,0.1,0.5,0.2]$
	
	$x^{*}(e_{4})=[0,0.75,0.7,1,0.75,0.4,0.1,0,0.2]$
	
	Also, the objective values of the above local optimal solutions are computed as $f_{1}\left(x^{*}\left(e_{1}\right)\right)=-0.9$, $f_{1}\left(x^{*}\left(e_{2}\right)\right)=-3.6$, $f_{1}\left(x^{*}\left(e_{3}\right)\right)=-1.3$ and $f_{1}\left(x^{*}\left(e_{4}\right)\right)=-3.3$. So, by Theorem \ref{thm-3}, the global minimum solution of the problem is $x^{*}\left(e_{2}\right)$ with an optimal value of $f_{1}\left(x^{*}\left(e_{2}\right)\right)=-3.6$.
\end{example}

\begin{example}[\textbf{Support function of a set}]\label{ex-5}
	Consider the objective function as follows:
	$$
	f_{2}(x)=\sup \left\{x^{T} y: y \in \mathbb{C}\right\}
	$$
	with domain $\text{\textbf{dom}~}f(x)=\left\{x \in  \mathbb{R}^{n}: \sup \left\{x^{T} y: y \in \mathbb{C}\right\}<\infty\right\}$, where $\mathbb{C} \subseteq  \mathbb{R}^{n}$ and $\mathbb{C} \neq \varnothing$. Many theoretical and applied aspects of this function can be found in \cite{1,2}. If $x$ is an arbitrary fixed point in $ \mathbb{R}^{n}$ and $\mathbb{C}$ is a non-empty bounded polyhedral set (polytope), then $\sup \left\{x^{T} y: y \in \mathbb{C}\right\}$ becomes equivalent to the following traditional linear programming problem:
	$$
	\begin{array}{ll}
		\max & z(y)=x^{T} y \\
		& y \in \mathbb{C}
	\end{array}
	$$
	where $x_{j}$'s and $y_{j}$'s ($j \in \mathscr{J}$) are the cost coefficients and the decision variables, respectively. Now, suppose that the cost vector $x$ is perturbed along the $j$'th unit vector $u_{j}$ (a vector having zero components, except for a 1 in the $j$'th position), that is, $x$ is replaced by $x+\lambda u_{j}$ where $\lambda \geq 0$ . In other words, we increase $x_{j}$ to $x_{j}+\lambda$ and hold other cost coefficients fixed. However, since the parametric analysis on the cost vector in a linear programming problem always produces a piecewise-linear and increasing function $z(\lambda)$ \cite{1}, it is concluded that $f_{2}$ is non-decreasing in variable $x_{j}$, $\forall j \in \mathscr{J}$. Consequently, we have $\mathscr{J}^{+}=\mathscr{J}$ and $\mathscr{J}^{-}=\varnothing$. Therefore, by considering the problem stated in Examples \ref{ex-1},\ref{ex-2} and \ref{ex-3}, the local optimal solutions $x^{*}\left(e_{1}\right)$, $x^{*}\left(e_{2}\right)$, $x^{*}\left(e_{3}\right)$ and $x^{*}\left(e_{4}\right)$ are obtained by (\ref{eq-9}) as follows:
	
	$x^{*}(e_{1})=[0,0.75,0.1,0,0.75,0.4,0.1,0.8,0.2]$
	
	$x^{*}(e_{2})=x^{*}\left(e_{4}\right)=[0,0.75,0.1,0,0.75,0.4,0.1,0,0.2]$
	
	$x^{*}(e_{3})=[0,0.75,0.1,0,0.75,0.4,0.1,0.5,0.2]$
	
	Also, for calculating the objective function values, let $\mathbb{C}$ be the polytope defined as the set of all the points $x \in  \mathbb{R}^{9}$ such that $\sum_{j=1}^{9} x_{j} \leq 1$ and $x_{j} \geq 0$, $\forall j \in\{1, \ldots, 9\}$. So, $f_{2}(x^{*}(e_{1}))$ is calculated as the optimal value of the following linear programming problem:
	$$
	\begin{array}{ll}
		\max & x^{*}\left(e_{1}\right)^{T} y \\
		& y \in \mathbb{C}
	\end{array}
	$$
	So, $f_{2}\left(x^{*}\left(e_{1}\right)\right)=0.8$ and in the same way, we get $f_{2}\left(x^{*}\left(e_{2}\right)\right)=0.75$, $f_{2}\left(x^{*}\left(e_{3}\right)\right)=0.75$ and $f_{2}\left(x^{*}\left(e_{4}\right)\right)=0.75$ . Hence, $x^{*}(e_{2})$, $x^{*}(e_{3})$ and $x^{*}(e_{4})$ are the global optimal solutions with the optimal value of 0.75.
\end{example}

\begin{example}[\textbf{perspective function}]\label{ex-6}
	Consider the following function as the objective function of the previous problem:
	$$
	f_{3}(x)=\frac{\left|x_{1}\right|^{p}+\cdots+\left|x_{8}\right|^{p}}{x_{9}^{p-1}}
	$$
	It is to be noted that since $I_{9}=[0.2,1]$ (from Table 3), then $0<0.2 \leq x_{9}$ for any $x \in S\left(A^{+}, A^{-}, b\right)$. So, $f_{3}$ is indeed the perspective function of $\left|x_{1}\right|^{p}+\cdots+\left|x_{8}\right|^{p}, p \geq 1$ \cite{2}. In this case, $\mathscr{J}^{+}=\{1, \ldots, 8\}$ and $\mathscr{J}^{-}=\{9\}$. Based on Definition \ref{def-8}, the local optimal solutions $x^{*}\left(e_{1}\right)$, $x^{*}\left(e_{2}\right)$, $x^{*}\left(e_{3}\right)$ and $x^{*}\left(e_{4}\right)$ are obtained for $p=3$ as follows:
	
	$x^{*}\left(e_{1}\right)=[0,0.75,0.1,0,0.75,0.4,0.1,0.8,1]$
	
	$x^{*}\left(e_{2}\right)=[0,0.75,0.1,0,0.75,0.4,0.1,0,0.5]$
	
	$x^{*}\left(e_{3}\right)=[0,0.75,0.1,0,0.75,0.4,0.1,0.5,0.2]$
	
	$x^{*}\left(e_{4}\right)=[0,0.75,0.1,0,0.75,0.4,0.1,0,0.2]$
	
	with objective values $f_{3}(x^{*}(e_{1}))=1.4218, f_{3}(x^{*}(e_{2}))=3.639$, $f_{3}(x^{*}(e_{3}))=25.8688$ and $f_{3}(x^{*}(e_{4}))=22.7438$. So, the global minimum solution of the problem is $x^{*}(e_{1})$ with an optimal value of $f_{3}(x^{*}(e_{1}))=1.4218$.
\end{example}

\begin{example}\label{ex-7}
	The following functions are among a wide range of functions that appear in many applications \cite{1,2}, and their optimal solutions are provided by Theorem 3:
	
	\begin{itemize}
		\item \textbf{Max function}. $f_{4}(x)=\max \left\{x_{1}, x_{2}, \ldots, x_{9}\right\}$.
		
		\item \textbf{Geometric mean}. $f_{5}(x)=\left(\prod_{j=1}^{9} x_{j}\right)^{1 / 9}$.
		
		\item \textbf{Log-sum-exp}. $f_{6}(x)=\operatorname{Ln}\left(e^{x_{1}}+e^{x_{2}}+\cdots+e^{x_{9}}\right)$. This function can be interpreted as a differentiable approximation of the max function, since $\max _{j=1}^{n}\left\{x_{j}\right\} \leq \operatorname{Ln}\left(\sum_{j=1}^{n} e^{x_{j}}\right) \leq \max _{j=1}^{n}\left\{x_{j}\right\}+\operatorname{Ln}(n)$ for all $x$.
		
		\item \textbf{Norms}. $f_{7}(x)=\left(\left|x_{1}\right|^{8}+\left|x_{2}\right|^{8}+\cdots+\left|x_{9}\right|^{8}\right)^{1 / 8}$.
		
		\item \textbf{Norm of a matrix}. $f_{8}(x)=\|x\|_{F}=\left(\sum_{j=1}^{9} x_{j}^{2}\right)^{1 / 2}$, where
		$
		x=\left[\begin{array}{lll}
			x_{1} & x_{2} & x_{3} \\
			x_{4} & x_{5} & x_{6} \\
			x_{7} & x_{8} & x_{9}
		\end{array}\right]
		$
		\item \textbf{Sum of four largest components}. $f_{9}(x)=\sum_{j=1}^{4} x_{[j]}$, where $x_{[j]}$ denotes the $j$'th largest component of $x$, i.e., $x_{[1]} \geq x_{[2]} \geq \cdots \geq x_{[9]}$ are the components of $x$ sorted in non-increasing order.
		
		\item \textbf{Maximum eigenvalue of a symmetric matrix}. The function $f_{10}(x)=\lambda_{\max }(x)$, where
		$
		x=\left[\begin{array}{lll}
			x_{6} & x_{1} & x_{2} \\
			x_{1} & x_{8} & x_{3} \\
			x_{2} & x_{3} & x_{9}
		\end{array}\right]
		$
		
		This function can be equivalently expressed as $f_{10}(x)=\sup \left\{u^{T} x u:\|u\|_{2}=1\right\}$, where $\|\cdot\|_{2}$ denotes the Euclidean norm.
		
		\item \textbf{Sum of logarithms}. $f_{11}(x)=\sum_{j=1}^{9} \operatorname{Ln}\left(\alpha_{j}+x_{j}\right)$, where $\alpha_{j}>0$. This function arises in information theory, in allocating power to a set of $n$ communication channels. The variable $x_{j}$ represents the transmitter power allocated to the $j$'th channel, and $\operatorname{Ln}\left(\alpha_{j}+x_{j}\right)$ gives the capacity or communication rate of the channel, so the problem is to allocate a total power of one to the channels, in order to maximize the total communication rate.
	\end{itemize}
	
	For all the above functions, we have $\mathscr{J}^{+}=\mathscr{J}$ and $\mathscr{J}^{-}=\varnothing$, and therefore the local optimal solutions $x^{*}\left(e_{1}\right)$, $x^{*}\left(e_{2}\right)$, $x^{*}\left(e_{3}\right)$ and $x^{*}\left(e_{4}\right)$ are obtained as in Example 5 , where $x^{*}\left(e_{2}\right)=x^{*}\left(e_{4}\right)$. Table \ref{t-6} shows the objective values of the solutions $x^{*}\left(e_{1}\right)$, $x^{*}\left(e_{2}\right)$ and $x^{*}\left(e_{3}\right)$ for $f_{i}(x)$, $i \in\{3,4, \ldots, 11\}$, where the results that are marked with "*" indicate the global optimal function values.
	\begin{table}[!ht]
		\centering
		\caption{The objective values of $f_{i}\left(x^{*}\left(e_{k}\right)\right)$ for $i \in\{3,4, \ldots, 11\}$ and $k \in\{1, \ldots, 4\}$.}
		\label{t-6}
		\scalebox{0.99}{
			\begin{tabular}{|c|c|c|c|}
				\hline
				& $f_{i}\left(x^{*}\left(e_{1}\right)\right)$ & $f_{i}\left(x^{*}\left(e_{2}\right)\right)$ & $f_{i}\left(x^{*}\left(e_{3}\right)\right)$ \\
				\hline
				$f_{4}$ & 0.8 & $0.75^{*}$ & $0.75^{*}$ \\
				\hline
				$f_{5}$ & $0^{*}$ & $0^{*}$ & $0^{*}$ \\
				\hline
				$f_{6}$ & 2.594 & $2.498^{*}$ & 2.5499 \\
				\hline
				$f_{7}$ & 0.8827 & $0.8182^{*}$ & 0.8201 \\
				\hline
				$f_{8}$ & 1.4089 & $1.1597^{*}$ & 1.2629 \\
				\hline
				$f_{9}$ & 2.7 & $2.1^{*}$ & 2.4 \\
				\hline
				$f_{10}$ & 1.0728 & $1.0607^{*}$ & 1.0644 \\
				\hline
				$f_{11}$ & 21.0238 & $20.9468^{*}$ & 20.9956 \\
				\hline
		\end{tabular}}
	\end{table}
\end{example}

% ==================================== conclusion ====================================
\section*{Conclusion}\label{sec_con}
This paper develops bipolar fuzzy relational equations for the most general case in which fuzzy compositions are defined by arbitrary continuous t-norms. An algorithm was proposed for finding a global optimal solution to a wide class of nonlinear optimization problems constrained by continuous bipolar fuzzy relation equations. The analytical properties of the continuous bipolar FREs were investigated, and the feasible solution set of the problem was completely characterized by a finite number of compact sets that may not be connected sets. Two necessary feasibility conditions and a necessary and sufficient condition were also presented to determine the feasibility of the problem. In addition, five simplification rules were proposed to speed up solution finding. Furthermore, as mentioned earlier, the proposed method provides some mathematical background that can be used to develop improved methods for solving optimization problems that are more likely to have different types of objective functions rather than decreasing (non-increasing) functions. Despite this, since the complexity of solving the problem increases exponentially with the problem size, it is always a challenge to find other simplification rules or methods that can reach the optimal solution.

%\newpage

% ==================================== sectionA_appendix ====================================

\section*{Appendix A}\label{sec_appA}
\setcounter{table}{0} \renewcommand{\thetable}{A\arabic{table}}

\begin{table}[!ht]
	\centering
	\caption{Summary of some commonly used continuous t-norms.}
	\label{tbl_a1}
	% \resizebox{\columnwidth}{!}{%
		\scalebox{0.9}{
			\begin{tblr}{c|c|c}
				\textbf{t-norm}  & \textbf{Function}                                                                                                                     & \textbf{Parameter} \\ \hline \hline
				Minimum & $T_M (x,y)=\min \{x,y\}$ & - \\ \hline
				Product          & $T_P(x, y)=xy$                                                                                                                        & -                  \\ \hline
				Einstein Product & $T_{EP}(x, y)=\frac{x y}{2-(x+y-x y)}$                                                                                                & -                  \\ \hline
				Lukasiewicz      & $T_{L}(x, y)=max\{0,x+y-1\}$                                                                                                          & -                  \\ \hline
				Frank            & $T_F^S(x, y)=\log _s\left(1+\frac{\left(s^x-1\right)\left(s^y-1\right)}{s-1}\right)$                                                  & $s>0, s\neq 1$     \\ \hline
				Yager            & $T_Y^P(x, y)=max\{0,1-[(1-x)^p+(1-y)^p]^{1/p}\}$                                                                                      & $P>0$              \\ \hline
				Hamacher         & $
				T_H^\alpha(x, y)=
				\begin{cases}
					0                                      & , \quad \alpha=x=y=0        \\
					\frac{x y}{\alpha+(1-\alpha)(x+y-x y)} & , \quad \text { otherwise }
				\end{cases}
				$                & $\alpha\geq 0$                                                                                                                                             \\ \hline
				Dombi            & $T_D^\lambda(x, y)=
				\begin{cases}
					0                                                                                                            & , \quad x=0 \text { or } y=0 \\
					\frac{1}{1+\left[\left(\frac{1-x}{x}\right)^\lambda+\left(\frac{1-y}{y}\right)^\lambda\right]^{1 / \lambda}} & , \quad \text {otherwise}
				\end{cases}
				$                & $\lambda>0$                                                                                                                                                \\ \hline
				Schweizer-Sklar  & $T_{S S}^p(x, y)=\left(\max \left\{0, x^p+y^p-1\right\}\right)^{1 / p}$                                                               & $p \neq 0$              \\ \hline
				Sugeno-Weber     & $T_{SW}^{\lambda}(x,y)=max\{0,(x+y-1+\lambda xy)/(1+ \lambda)\}$                                                                      & $\lambda>-1$       \\ \hline
				Aczel-Alsina     & $T_{A A}^\lambda(x, y)=\exp \left(-\left[(-\operatorname{Ln}(x))^\lambda+(-\operatorname{Ln}(y))^\lambda\right]^{1 / \lambda}\right)$ & $\lambda>0$ \\ \hline
				Dubois-Prade & $T^{\gamma}_{DP}=xy/\max\{x,y,\gamma\}$ & $0 \leq \gamma \leq 1$ \\ \hline
				Mayor-Torrence & $
				T_{MT}^{\lambda}(x, y)=
				\begin{cases}
					\max \{0,x+y+\lambda\}                                      & , ~ 0<\lambda \leq 1,~(x,y) \in {{[0,\lambda]}^{2}}        \\
					\min \{x,y\} & , ~ \text { otherwise }
				\end{cases}
				$ & $0 \leq \lambda \leq 1$
			\end{tblr}%
		}
	\end{table}
	
	\begin{table}[!ht]
		\centering
		\caption{Values $l$ and $u$ discussed in Lemma \ref{lm-1} for the t-norms stated in Table \ref{tbl_a1}. The values $l$ and $u$ are presented for the cases that $a \geq b>0$ and $a>b \geq 0$, respectively; because, for all the cases, it is easily verified that $l=0$ if $b=0$, and $u=1$ if $a=b$.}
		\label{tbl_a2}
			\scalebox{0.9}{
				\begin{tblr}{c|c|c}
					t-norm & $l(a \geq b>0)$ & $u(a>b \geq 0)$ \\
					\hline \hline
					$T_{M}$ & $b$
					& $b$ \\
					\hline
					$T_{P}$ & $b / a$ & $b / a$ \\
					\hline
					${T_{E P}}$ & $((2-a) b) /(a+b-a b)$ & $((2-a) b) /(a+b-a b)$ \\
					\hline
					$T_{L}$ & $1+b-a$ & $1+b-a$ \\
					\hline
					$T_{F}^{s}$ & $\log _{s}\left(1+\left[(s^{b}-1)(s-1) /\left(s^{a}-1\right)\right]\right)$ & $\log _{s}\left(1+\left[(s^{b}-1)(s-1) /\left(s^{a}-1\right)\right]\right)$ \\
					\hline
					$T_{Y}^{p}$ & $1-\left((1-b)^{p}-(1-a)^{p}\right)^{1 / p}$ & $1-\left((1-b)^{p}-(1-a)^{p}\right)^{1 / p}$ \\
					\hline
					$T_{H}^{\alpha}$ & $([\alpha+(1-\alpha) a] b) /(a-(1-\alpha)(1-a) b)$ & $([\alpha+(1-\alpha) a] b) /(a-(1-\alpha)(1-a) b)$ \\
					\hline
					$T_{D}^{\lambda}$ & $\left(1+\left[((1-b) / b)^{\lambda}-((1-a) / a)^{\lambda}\right]^{1 / \lambda}\right)^{-1}$ & $\begin{cases}\left(1+\left[((1-b) / b)^{\lambda}-((1-a) / a)^{\lambda}\right]^{1 / \lambda}\right)^{-1} & , b>0 \\ 0 & , b=0\end{cases}$ \\
					\hline
					${T_{S S}^{p}}$ & $\left(1+b^{p}-a^{p}\right)^{1 / p}$ & $\begin{cases}\left(1+b^{p}-a^{p}\right)^{1 / p} & , b>0 \\ 0 & , b=0, p<0\end{cases}$ \\
					\hline
					$T_{S W}^{\lambda}$ & $((1+\lambda) b+1-a) /(1+\lambda a)$ & $((1+\lambda) b+1-a) /(1+\lambda a)$ \\
					\hline
					$T_{A A}^{\lambda}$ & $\exp \left(-\left[(-\operatorname{Ln}(b))^{\lambda}-(-\operatorname{Ln}(a))^{\lambda}\right]^{1 / \lambda}\right)$ & $\begin{cases}\exp \left(-\left[(-\operatorname{Ln}(b))^{2}-(-\operatorname{Ln}(a))^{\lambda}\right]^{1 / \lambda}\right) & , b>0 \\ 0 & , b=0\end{cases}$ \\
					\hline
					$T_{D P}^{\gamma}$ & $\begin{cases}\max \{b, \gamma\} & , a=b \\ \gamma b / a & , b<a<\gamma \\ b & , a>b, a \geq \gamma\end{cases}$ & $\begin{cases}\gamma b / a & , b<a<\gamma \\ b & , a>b, a \geq \gamma\end{cases}$ \\
					\hline
					$T_{M T}^{\lambda}$ & $\begin{cases}\lambda & , a=b, a \leq \lambda \\ b & , a=b, a>\lambda \\ b+\lambda-a & , a>b, a \leq \lambda \\ b & , a>b, a>\lambda\end{cases}$ & $\begin{cases}b+\lambda-a & , a>b, a \leq \lambda \\ b & , a>b, a>\lambda\end{cases}$ 
			\end{tblr}}
		\end{table}

		%============================AppendixB
		\section*{Appendix B}\label{sec_appB}
		All the proofs discussed in this section are based on the fact that under some assumptions, any solution $x \in S_{\{i\}}$ $\left(A^{+}, A^{-}, b\right)$ also satisfies the $i$'th equation of Problem (\ref{eq_1}), i.e., $S_{\left\{i_{0}\right\}}\left(A^{+}, A^{-}, b\right) \subseteq S\left(A^{+}, A^{-}, b\right)$, and therefore $S\left(A^{+}, A^{-}, b\right)=S_{\left\{i_{0}\right\}}\left(A^{+}, A^{-}, b\right)$. Consequently, if this is the case, the $i$'th equation is indeed a redundant (irrelevant) constraint (i.e., a constraint that does not affect the feasible region), and therefore it can be deleted from further consideration.
		
		\renewcommand{\thetheorem}{B1} % or \Alph{section}.\arabic{theorem}

		\begin{lemma}\label{lm-b1}
			Suppose that $S\left(A^{+}, A^{-}, b\right) \neq \varnothing$ and $i_{0} \in \mathscr{I}$. If $x \in S_{\left\{i_{0}\right\}}\left(A^{+}, A^{-}, b\right)$, then $x_{j} \in I_{j}$, $\forall j \in \mathscr{J}$. Particularly, $x_{j} \in I_{i_{0} j}$, $\forall j \in \mathscr{J}$.
		\end{lemma}
		
		\begin{proof}
			From Theorem \ref{thm-1}, $S_{\left\{i_{0}\right\}}\left(A^{+}, A^{-}, b\right)=\bigcup_{e^{\prime} \in E^{\prime}} S\left(e^{\prime}\right)$, where $E^{\prime}$ is the set of all the restrictions of admissible functions $e \in E$ to $\mathscr{I}-\left\{i_{0}\right\}$. So, for each $x \in S_{\left\{i_{0}\right\}}\left(A^{+}, A^{-}, b\right)$, there exists at least one $e_{0}^{\prime} \in E^{\prime}$ such that $x \in S\left(e_{0}^{\prime}\right)$. Now, from (\ref{eq-5}), for each $j \in \mathscr{J}$ we have either $x_{j} \in \bigcap_{i \in \mathscr{I}_{j}\left(e_{0}^{\prime}\right)} S_{i j}^{\prime}$ or $x_{j} \in I_{j}$. But, in the former case, since $S_{i j}^{\prime}=S_{i j} \cap I_{j}$ (Definition \ref{def-4}), we have again $x_{j} \in \bigcap_{i \in \mathscr{I}_{j}\left(e_{0}^{\prime}\right)} S_{i j}^{\prime} \subseteq I_{j}$. Consequently, $x_{j} \in I_{j}$ (and therefore $x_{j} \in I_{i_{0} j}$ from Definition \ref{def-4}), $\forall j \in \mathscr{J}$.
		\end{proof}
		
		\begin{proof}[The Proof of Rule \ref{rule-1}]
	Let $x \in S_{\left\{i_{0}\right\}}\left(A^{+}, A^{-}, b\right)$. We shall show that $x$ also satisfies the $i_{0}$'th equation. From Lemma \ref{lm-b1}, $x_{j} \in I_{j}$, $\forall j \in \mathscr{J}$(*1). Also, since $S\left(A^{+}, A^{-}, b\right) \neq \varnothing$, there exists at least one $j_{0} \in \mathscr{J}$ such that $S_{i_{0} j_{0}}^{\prime} \neq \varnothing$ (Lemma \ref{lm-3} (b)). On the other hand, since $b_{i_{0}}=0$, then we have $S_{i_{0} j_{0}}^{\prime}=I_{i_{0} j_{0}}$ which together with $I_{j_{0}} \subseteq I_{i_{0} j_{0}}$ (Definition \ref{def-4}) and (*1) imply $x_{j_{0}} \in S_{i_{0} j_{0}}^{\prime}$ (*2). Now, the result follows from (*1),(*2) and Corollary \ref{corl-2}.
\end{proof}

\begin{proof}[The Proof of Rule \ref{rule-2}]
	For any $x \in S\left(A^{+}, A^{-}, b\right)$, Lemma \ref{lm-4} implies $x_{j_{0}} \in I_{j_{0}}=\{k\}$, i.e., $x_{j_{0}}=k$. Now, from Lemma \ref{lm-b1}, if $x \in S_{\left\{i_{0}\right\}}\left(A^{+}, A^{-}, b\right)$, then $x_{j} \in I_{i_{0} j}$, $\forall j \in \mathscr{J}$(*1). But, since $k \in S_{i_{0} j_{0}}^{\prime}$, then we have $S_{i_{0} j_{0}}^{\prime} \neq \varnothing$, that together with $S_{i_{0} j_{0}}^{\prime}=S_{i_{0} j_{0}} \cap I_{j_{0}}$ and $I_{j_{0}}=\{k\}$ imply $S_{i_{0} j_{0}}^{\prime}=\{k\}$. Consequently, $x_{j_{0}} \in S_{i_{0} j_{0}}^{\prime}$ (*2). Hence, from (*1),(*2) and Corollary \ref{corl-2}, it follows that $x$ satisfies the $i_{0}$'th equation.
\end{proof}

\begin{proof}[The Proof of Rule \ref{rule-3}]
	If $x \in S_{\left\{i_{0}\right\}}\left(A^{+}, A^{-}, b\right)$, then by Lemma \ref{lm-b1}, we conclude that $x_{j} \in I_{i_{0} j}$, $\forall j \in \mathscr{J}$ (*1). Also, from Lemma \ref{lm-4}, there exists at least one $j_{i} \in \mathscr{J}$ such that $x_{j_{i}} \in S_{i j_i}^{\prime}$. But, since $S_{i j}^{\prime} \subseteq S_{i_{0} j}^{\prime}$, then we have $x_{j_{i}} \in S_{i_{0} j_{i}}^{\prime}$ (*2). So, (*1), (*2) and Corollary \ref{corl-2} imply that $x$ satisfies the $i_{0}$'th equation.
\end{proof}

\begin{proof}[The Proof of Rule \ref{rule-4}]
	If $x \in S\left(A^{+}, A^{-}, b\right)$, then from Lemma \ref{lm-4}, there exists at least one $j \in \mathscr{J}$ such that $x_{j} \in S_{i_{0} j}^{\prime}$. But, since $\mathscr{J}_{i_{0}}=\left\{j_{0}\right\}$, we have necessarily $x_{j_{0}} \in S_{i_{0} j_{0}}^{\prime}=\{k\}$, that means $x_{j_{0}}=k$. Moreover, suppose that $i \in \mathscr{I}$ and $k \in S_{i j_{0}}^{\prime}$. If $x \in S_{\{i\}}\left(A^{+}, A^{-}, b\right)$, then $x$ satisfies the $i_{0}$'th equation, and therefore we must have $x_{j_{0}}=k \in S_{i j_{0}}^{\prime}$. Also, from Lemma \ref{lm-b1}, $x_{j} \in I_{i j}$, $\forall j \in \mathscr{J}$. Hence, Corollary \ref{corl-2} implies that $x$ also satisfies the $i$'th equation.
\end{proof}

\begin{proof}[The Proof of Rule \ref{rule-5}]
	Let $x \in S_{\left\{i_{0}\right\}}\left(A^{+}, A^{-}, b\right)$. From Lemma \ref{lm-b1}, $x_{j_{0}} \in I_{j_{0}}$ and $x_{j} \in I_{i_{0} j}$, $\forall j \in \mathscr{J}$. But, since $S_{i_{0} j_{0}}^{\prime}=I_{j_{0}}$, then $x_{j_{0}} \in I_{j_{0}}=S_{i_{0} j_{0}}^{\prime}$, and therefore from Corollary \ref{corl-2} it follows that $x$ satisfies the $i_{0}$'th equation.
\end{proof}

% ==================================== refrences ====================================

\section*{Acknowledgement}
The authors wish to express their appreciation for several excellent suggestions
for improvements in this paper made by the referees.\\

%\newpage
%\bibliographystyle{ieeetr} 
\bibliography{refrences}
%\printbibliography

\end{document}